\documentclass[a4paper, 11pt, oneside, onecolumn,sort&compress]{elsarticle}
\usepackage{amssymb}
\usepackage{amsthm}
\usepackage{mathrsfs}
\usepackage{amsfonts}
\usepackage[title]{appendix}
\usepackage{syntonly}
\usepackage{xcolor} 
\usepackage{lineno}

\newtheorem{theorem}{Theorem}[section]
\newtheorem{lemma}{Lemma}[section]

\newdefinition{remark}{Remark}[section]
\newdefinition{example}{Example}[section]
\newdefinition{proposition}{Proposition}[section]
\newdefinition{definition}{Definition}[section]

\usepackage{amssymb,amsmath}

\keywordtitle{\normalfont\bfseries Keywords}

\numberwithin{equation}{section}
\makeatletter
\def\ps@pprintTitle{%
  \let\@oddhead\@empty
  \let\@evenhead\@empty
  \let\@oddfoot\@empty
  \let\@evenfoot\@oddfoot}
\makeatother

\begin{document}
	

\begin{frontmatter}

\title{\textbf{Turing instability for periodic reaction-diffusion systems}
	\tnoteref{label1}}

\author[mainaddress] {Yequan Liang\corref{cor1}}
\ead{yequanl@mun.ca}

\author[mainaddress] {Xiao-Qiang Zhao\corref{cor1}}
\ead{zhao@mun.ca}

\address[mainaddress]{Department of Mathematics and Statistics, Memorial University of Newfoundland, \\[0.5em]
St. John's, NL A1C 5S7, Canada}

\begin{abstract}
This paper is devoted to the study of Turing instability of spatially
homogeneous periodic solutions for time-periodic reaction--diffusion systems with Neumann boundary conditions. We first identify three classes of such  periodic systems for which the Turing instability does not 
occur.   Then we analyze a class of periodic Riccati equations and give  a complete classification of their global dynamics.  This result enables us to establish a new criterion for Turing instability in periodic  
reaction-diffusion systems of predator-prey type.  A further comparison with  the autonomous case shows that our Turing instability criterion is equivalent to  the classical one for spatially
homogeneous equilibria, while it is complementary to a criterion for the Hopf-Turing instability.
\end{abstract}

\medskip
\end{frontmatter}


\section{Introduction }

Since Turing's pioneering work \cite{Turing-PTRSLB-1952}, there have
been extensive investigations of diffusion-driven instability
(i.e., Turing instability) in reaction--diffusion systems, see, e.g., 
\cite{Casten-Holland-SIAP-1977, Sherratt-JMB-1995, 
	Satnoianu-Menzinger-Maini-JMB-2000, Ni-Tang-TAMS-2005, Ricard-Mischler-JNS-2009, Shi-Xie-Little-JAAC-2011,
	Yi-JDE-2021} and the references therein. This phenomenon occurs when a spatially homogeneous state is stable for the kinetic system but becomes unstable after diffusion is introduced. It provides a mechanism for the emergence of spatial patterns in biological and chemical systems
\cite{Murray-book-2003}. Thus, diffusion can promote rather than merely suppress spatial heterogeneity.

Much of the existing literature concerns autonomous reaction--diffusion systems. For spatially homogeneous equilibria, the classical theory reduces
the stability problem to an algebraic eigenvalue problem
\cite{Murray-book-2003}. Casten and Holland \cite{Casten-Holland-SIAP-1977} rigorously reduced the linear stability of spatially homogeneous equilibria
under Neumann boundary conditions to a family of finite-dimensional
eigenvalue problems associated with the Neumann Laplacian.
Extensions to multicomponent systems and cross-diffusion were studied in
\cite{Satnoianu-Menzinger-Maini-JMB-2000,Shi-Xie-Little-JAAC-2011}.
Beyond equilibria, Maginu \cite{Maginu-JDE-1979} obtained a
sufficient condition for the instability of spatially homogeneous
periodic solutions in terms of the variation of their period under
an auxiliary perturbation. Building on Maginu's approach, Yi \cite{Yi-JDE-2021} derived sufficient conditions for the Turing instability of stable periodic
solutions arising from a Hopf bifurcation, with extensions to reaction--diffusion systems with cross-diffusion and related patch models.
Turing instability near Hopf bifurcation has also been investigated in a number of general and model-specific settings \cite{Dilao-AMC-2005,Yi-Wei-Shi-NARWA-2008,Ricard-Mischler-JNS-2009,
	Li-Jiang-Shi-IMA-2013,Wu-Zhou-Shao-Chen-PhysA-2017, Wu-Yi-MA-2026}.

However, Turing instability in time-periodic reaction--diffusion systems 
has received comparatively less attention. Timm and Okubo \cite{Timm-Okubo-JMB-1992} investigated a plankton predator--prey model with periodically varying diffusivity.
Dancer \cite{Dancer-DIE-1994} constructed a periodic predator-prey reaction--diffusion system of Lotka-Volterra type for which the Turing instability occurs. Sherratt \cite{Sherratt-JMB-1995,Sherratt-EJAM-1995} considered diffusion and reaction parameters that alternate periodically
between constant values and derived instability conditions
in regimes of slow or fast temporal oscillations. Related effects of time-periodic diffusivities were further
studied in \cite{Bhattacharyya-Mukhopadhyay-Bandyopadhyay-IJMMS-2003, Mukhopadhyay-Bhattacharyya-BMB-2006}. 
Wang and Lutscher \cite{Wang-Lutscher-JMB-2019} obtained
diffusion-driven instability conditions for an impulsive
predator--prey system in a two-season environment. Hetzer, Madzvamuse, and Shen
\cite{Hetzer-Madzvamuse-Shen-DCDS-2012}
characterized Turing instability on evolving domains through
Lyapunov exponents of the associated evolution family.
More recently, Singh et al.
\cite{Singh-Bhargava-Dubey-ND-2025} studied Turing instability
in a seasonally forced predator--prey model with time-dependent
cross-diffusion and periodic reaction rates. In a more general nonautonomous framework, Van Gorder \cite{VanGorder-PRSA-2020} derived sufficient instability conditions for two-component reaction--diffusion systems using comparison principles, without specializing the analysis to
time-periodic systems. Thus, more refined and verifiable criteria for Turing instability in 
time-periodic reaction--diffusion systems need to be further explored. 

The purpose of this paper is to study the  Turing instability for general  time-periodic reaction--diffusion systems. The major difficulty lies in the fact that one cannot explicitly solve the Floquet multipliers from the linearized periodic ODE system (even for two-dimensional one). This is quite different from the Turing instability of a periodic orbit for an autonomous system \cite{Ricard-Mischler-JNS-2009, Yi-JDE-2021}, where one of the Floquet multipliers for the linearized periodic ODE system is always one, but this is not true, in general, for the linearized periodic ODE system from a periodic reaction-diffusion system. Following the modal approach in \cite{Casten-Holland-SIAP-1977,Dancer-DIE-1994}, the linearized
stability for a spatially homogeneous periodic solution on
a bounded domain with homogeneous Neumann boundary conditions is equivalent to that for 
an infinite family of linear periodic ODE systems indexed
by the Neumann eigenvalues of $-\Delta$.
It turns out that the Turing instability is reduced to determining whether one of
the nonzero spatial modes has a Floquet multiplier outside the unit circle.

In this paper, we first show that Turing instability does not occur
when the periodic reaction-diffusion system admits a comparison principle, when all diffusion
coefficients are sufficiently large, and  when all diffusion rates are
equal, respectively. To study the occurrence of Turing instability, we turn to
periodic reaction-diffusion systems of predator--prey type. By analyzing a class of
scalar periodic Riccati equations, we obtain a complete classification
of their global dynamics, which is of independent interest, and identify a critical modal parameter governing the existence and multiplicity of positive periodic solutions.
For modal parameters at or above the critical value, a distinguished
positive periodic solution yields an exact integral representation
of the dominant Floquet multiplier. We further prove that the resulting growth
integral is strictly concave above the critical value and that the unstable modal set is either empty or a single bounded interval. Based on this analysis, we derive a sufficient criterion involving positive periodic supersolutions, which avoids the need to determine the periodic solution of the Riccati equation explicitly. These modal results yield Turing instability on a fixed domain whenever an unstable positive parameter is a nonzero Neumann
eigenvalue, and on a suitably rescaled domain otherwise. We further find that  our criterion  coincides with 
the classical one for the Turing instability in  two dimensional autonomous systems, and  complements a criterion obtained in  \cite{Yi-JDE-2021} for the Hopf-Turing instability on periodic orbits of autonomous reaction-diffusion systems.

The rest of this paper is organized as follows.
In the next section, we introduces the stability
framework and the modal reduction. In Section~\ref{pthirdsec-Nonoccurrence-of-Turing-Instability}, we
present  the nonoccurrence results under three cases. In Section 4, we develop Turing instability criteria based on a class of periodic Riccati equations and compares them with those obtained for autonomous systems. 
In Section~\ref{pthirdsubsec-applications},  we apply these analytical results  to the autonomous Sporns--Seelig model and the periodic Segel--Jackson model, respectively.

\section{Preliminaries}\label{pthirdsec-preliminaries}

In this section, we introduce the notation and basic facts used throughout the paper. We first formulate the periodic reaction-diffusion system and the corresponding kinetic system, and then specify the stability notion used in the definition of Turing instability. Finally, using the Neumann eigenfunction expansion, we reduce the linearized reaction-diffusion problem along a spatially homogeneous periodic solution to a family of finite-dimensional periodic systems.

Let $\omega>0$ be given and 
$\Omega\subset \mathbb R^N$ be a bounded connected domain with smooth boundary. We consider the following $\omega$-periodic reaction-diffusion system subject to the homogeneous Neumann boundary condition:
\begin{equation}\label{pthirdsys-nonlinear-PDE}
	\left\{
	\begin{aligned}
		\frac{\partial u}{\partial t} &= D(t)\Delta u+F(t,u),
		\qquad &&(t,x)\in (0,\infty)\times\Omega,\\
		\frac{\partial u}{\partial \nu} &=0,
		\qquad &&(t,x)\in (0,\infty)\times\partial\Omega.
	\end{aligned}
	\right.
\end{equation}
Here $u=(u_1,\ldots,u_m)$, $F(t,u)=(f_1(t,u),\ldots,f_m(t,u))$, $D(t)=\operatorname{diag}(d_1(t),\ldots,d_m(t))$ and  $\Delta u = (\Delta u_1, \ldots,\Delta u_m)$. We assume that each $d_i$ is positive, H\"older continuous, and $\omega$-periodic in $t$. The vector $\nu$ denotes the outward normal vector to $\partial \Omega$. In population biology, the Neumann boundary condition represents the fact that the population is confined in the region $\Omega$ and there is no migration across the boundary.

The kinetic system associated with system \eqref{pthirdsys-nonlinear-PDE} is the periodic ordinary differential equation:
\begin{equation}\label{pthirdsys-nonlinear-ODE}
	\frac{\mathrm d v}{\mathrm d t}=F(t,v).
\end{equation}
Clearly, every solution of system \eqref{pthirdsys-nonlinear-ODE} is naturally a spatially homogeneous solution of system \eqref{pthirdsys-nonlinear-PDE}. In what follows, we fix a periodic solution \(u^*(t)\) of the
kinetic system \eqref{pthirdsys-nonlinear-ODE} such that $u^*(t+\omega)=u^*(t)$ for all $t\in\mathbb R$, and denote by $\check u^*(t,x):=u^*(t)$ the corresponding spatially homogeneous solution of
\eqref{pthirdsys-nonlinear-PDE}. When the kinetic system is autonomous and \(u^*(t)\) is nonconstant, \(\omega\) is taken to be the least period of \(u^*(t)\). We impose the following regularity assumption on $F$.

\begin{itemize}
	\item[(A1)] The mapping $F \in C^1(\mathbb R\times\mathbb R^m;\mathbb R^m)$ is \(\omega\)-periodic in \(t\), and \(D_uF\) is locally Lipschitz continuous on \(\mathbb R\times\mathbb R^m\).
\end{itemize}

Let $\mathbb X:=C(\overline\Omega;\mathbb R^m)$, equipped with the supremum norm $\|\phi\|_{\mathbb X}=\max_{x\in\overline\Omega}|\phi(x)|$. We write $u(t,\phi)$ and $v(t,a)$ for the solutions of systems \eqref{pthirdsys-nonlinear-PDE} and \eqref{pthirdsys-nonlinear-ODE} with initial data $\phi$ and $a$ at $t=0$, respectively. For $a\in\mathbb R^m$, let $\check a(x):=a$. The subspace of spatially constant functions is invariant under
the reaction-diffusion dynamics, and the restricted dynamics coincide with those of the kinetic system.

Throughout this paper, the stability of $u^*(t)$ refers to the usual Lyapunov stability when the kinetic system~\eqref{pthirdsys-nonlinear-ODE} is nonautonomous, and to the orbital stability when it is autonomous. In the latter case, write $\Gamma:=\{u^*(t):0\le t<\omega\}$ and $\check\Gamma:=\{\check a:a\in\Gamma\}$. For equilibria, orbital stability coincides with Lyapunov stability. The same convention applies to $\check u^*(t)$ in $\mathbb X$. We use the Euclidean norm for the kinetic system and the supremum norm of $\mathbb X$ for the reaction-diffusion system.

\begin{definition}
	We say that a Turing instability occurs if $u^*(t)$ is stable for the kinetic system~\eqref{pthirdsys-nonlinear-ODE}, whereas $\check{u}^*(t)$ is unstable for the reaction-diffusion system~\eqref{pthirdsys-nonlinear-PDE} in $\mathbb X$.
\end{definition}

Since the spatially homogeneous subspace is invariant and its
restricted dynamics coincide with those of the kinetic system, such an instability necessarily involves spatially heterogeneous perturbations. For later use, we impose the following assumption:
\begin{itemize}
	\item[(A2)] $u^*(t)$ is stable for the kinetic system~\eqref{pthirdsys-nonlinear-ODE}.
\end{itemize}

We now linearize the kinetic and reaction-diffusion systems along $u^*(t)$. Let $A(t):=D_uF(t,u^*(t))$. Since both $F$ and $u^*$ are $\omega$-periodic in $t$, the matrix $A(t)$ is also $\omega$-periodic. Consider the linear periodic system
\begin{equation}\label{pthirdsys-linear-ODE-zero}
	\frac{\mathrm d y}{\mathrm d t}=A(t)y.
\end{equation}
Let $T_0(t,s)$ be its evolution operator and let $T_0:=T_0(\omega,0)$ be the associated monodromy matrix. In the autonomous nonconstant periodic-orbit case, differentiating the kinetic equation along \(u^*(t)\) gives $T_0\dot u^*(0)=\dot u^*(0)$. Thus \(1\) is the trivial Floquet multiplier associated with phase translation. Consequently, one cannot characterize kinetic stability in this case by the condition $\rho(T_0)<1$, where \(\rho(\cdot)\) denotes the spectral radius of an operator.

Similarly, the linearization of system \eqref{pthirdsys-nonlinear-PDE} at $\check u^*(t)$ is
\begin{equation}\label{pthirdsys-linear-PDE}
	\left\{
	\begin{aligned}
		\frac{\partial w}{\partial t}&=D(t)\Delta w+A(t)w,\\
		\frac{\partial w}{\partial\nu}&=0.
	\end{aligned}
	\right.
\end{equation}
Let $\mathcal T(t,s)$ be the evolution operator associated with system~\eqref{pthirdsys-linear-PDE}, and set $\mathcal T:=\mathcal T(\omega,0)$. The variation-of-constants formulation and Assumption~(A1)
yield well-defined local period maps $\mathcal P(\phi):=u(\omega,\phi)$ and $P_0(a):=v(\omega,a)$, which are continuously Fr\'echet differentiable in neighborhoods of $\check u^*(0)$ in $\mathbb X$ and $u^*(0)$ in $\mathbb R^m$,
respectively \cite[Corollary~3.1.24~(ii) and Chapters~6--7]
{Lunardi-book-1995}. Their derivatives at the reference initial data satisfy $D\mathcal P(\check u^*(0))=\mathcal T$ and $DP_0(u^*(0))=T_0$.

Accordingly, we focus below on the Floquet spectra of these linear periodic systems. Their connection to nonlinear Turing instability is supplied by the linearized instability principles used in Lemma~\ref{pthirdlem-instability-reduction}(iii).

We now reduce the linearized PDE problem to a family of finite-dimensional
periodic systems. Accordingly, in the sequel we focus on these finite-dimensional modal systems in our study of Turing instability. Let
\[
0=\lambda_0< \lambda_1\le \lambda_2\le \cdots,\qquad \lambda_n\to+\infty,
\]
be the eigenvalues of $-\Delta$ under the homogeneous Neumann boundary condition, counted with multiplicity. Let $\{\phi_n\}_{n=0}^\infty$ be the corresponding orthonormal basis of $L^2(\Omega)$:
\[
-\Delta\phi_n=\lambda_n\phi_n,\qquad
\frac{\partial\phi_n}{\partial\nu}=0.
\]
Since \(\Omega\) is connected, the eigenvalue \(\lambda_0=0\) is simple
and its eigenfunction is spatially constant. Hence \(n=0\) represents
homogeneous perturbations, whereas every mode \(n\ge 1\) is spatially
heterogeneous. For each \(n\ge0\), define $E_n:=\{\phi_n(x)y:y\in\mathbb R^m\}$.
When \(y\in\mathbb C^m\) is used below, \(\phi_n y\) is understood as an
element of the complexification of \(\mathbb X\). For any \(s\ge0\) and \(\xi\in\mathbb R^m\), the solution of the linearized equation with initial value \(w(s,x)=\phi_n(x)\xi\) is given by $w(t,x)=\phi_n(x)T_n(t,s)\xi$ for all  $t\ge s$, where \(T_n(t,s)\) is the evolution operator associated with
\begin{equation}\label{pthirdsys-linear-decomposition-discrete}
	\frac{\mathrm d y}{\mathrm d t}
	=
	\bigl(A(t)-\lambda_nD(t)\bigr)y.
\end{equation}
Consequently, \(E_n\) is invariant under the linearized evolution
\(\mathcal T(t,s)\). We write $T_n:=T_n(\omega,0)$. The case \(n=0\) corresponds exactly to the kinetic variational system \eqref{pthirdsys-linear-ODE-zero}. Since \(T_n(t,s)\) is the evolution operator of a finite-dimensional linear system, uniqueness of solutions implies that \(T_n\) is invertible for every \(n\ge0\). In particular, \(0\notin\sigma(T_n)\), where $\sigma(\cdot)$ denotes the spectrum of an operator.

Define 
\[
A^*:=\max_{t\in[0,\omega]}\|A(t)\|, \quad d_*:= \min_{t\in [0, \omega]}\{d_1(t), d_2(t),\dots , d_m(t)\}>0. 
\]
Here \(\|\cdot\|\) denotes the operator norm induced by the Euclidean norm on \(\mathbb C^m\). Unless otherwise stated, spectra of real operators are taken over \(\mathbb C\). We use the same notation for an operator and its complexification. The following lemma makes the modal reduction of \(\mathcal T\) precise.

\begin{lemma}\label{pthirdlem-instability-reduction}
Let  (A1) hold. Then the following statements are valid:
	\begin{enumerate}
		\item[(i)] The nonzero spectrum of the monodromy operator \(\mathcal T\) satisfies
		\[
		\sigma(\mathcal T)\setminus\{0\}
		=
		\bigcup_{n\ge0}\sigma(T_n).
		\]
		
		\item[(ii)] For every \(n\ge1\), there holds $\|T_n\|\le e^{(A^*-\lambda_nd_*)\omega}$. In particular, \(\rho(T_n)<1\) for all sufficiently large \(n\).
		\item[(iii)] If there exists \(n_0\ge1\) such that \(\rho(T_{n_0})>1\), then
		\(\check{u}^*(t)\) is unstable for system~\eqref{pthirdsys-nonlinear-PDE} in \(\mathbb X\).
	\end{enumerate}
\end{lemma}

\begin{proof}
We first prove (i). By parabolic regularity, $\mathcal T$ is compact on $\mathbb X$ \cite[Section~2 and p.~26]{Amann-book-1978}. Therefore every nonzero spectral value of $\mathcal T$ is an eigenvalue. For  \(\xi=(\xi_1,\ldots,\xi_m)\) in the complexification of \(\mathbb X\), define its \(n\)-th vector Fourier coefficient by
\[
	c_n(\xi):=
	\int_\Omega \xi(x)\phi_n(x)\,\mathrm dx
	=
	\left(
	\int_\Omega \xi_1(x)\phi_n(x)\,\mathrm dx,\ldots,
	\int_\Omega \xi_m(x)\phi_n(x)\,\mathrm dx
	\right)
	\in\mathbb C^m .
\]
This is well-defined because \(\mathbb X\hookrightarrow L^2(\Omega;\mathbb R^m)\), and hence the complexification of \(\mathbb X\) is continuously embedded in \(L^2(\Omega;\mathbb C^m)\). Let $w(t):=\mathcal T(t,0)\xi$. For \(t>0\), parabolic regularity makes \(w(t,\cdot)\) sufficiently regular. Hence, by integration by parts and the Neumann boundary condition, for every \(n\ge0\),
\[
	\begin{aligned}
		\frac{\mathrm d}{\mathrm dt}c_n(w(t))
		&=
		\int_\Omega
		\left[D(t)\Delta w(t,x)+A(t)w(t,x)\right]\phi_n(x)\,\mathrm dx  \\
		&=
		-\lambda_nD(t)c_n(w(t))+A(t)c_n(w(t))  \\
		&=
		\bigl(A(t)-\lambda_nD(t)\bigr)c_n(w(t)).
	\end{aligned}
\]
By continuity at $t=0$, it follows that $c_n(w(t))=T_n(t,0)c_n(\xi)$ for all $ t\ge0$. In particular, we have $c_n(\mathcal T\xi)=T_n c_n(\xi)$.

Let \(\mu\in\sigma(\mathcal T)\setminus\{0\}\). Choose a nonzero eigenfunction \(\xi\) in the complexification of \(\mathbb X\) such that $\mathcal T\xi=\mu\xi$. Since \(\{\phi_n\}_{n=0}^\infty\) is complete in \(L^2(\Omega;\mathbb C)\), componentwise completeness implies that, if \(c_n(\xi)=0\) for all \(n\ge0\), then \(\xi=0\) in \(L^2(\Omega;\mathbb C^m)\). Hence there exists \(n_*\ge0\) such that \(c_{n_*}(\xi)\ne0\). Therefore
\[
	T_{n_*}c_{n_*}(\xi)
	=
	c_{n_*}(\mathcal T\xi)
	=
	\mu c_{n_*}(\xi).
\]
Thus \(\mu\in\sigma(T_{n_*})\).
	
Conversely, fix \(n\ge0\) and suppose that \(\mu\in\sigma(T_n)\). Since
\(T_n\) is finite-dimensional, there exists \(y\in\mathbb C^m\setminus\{0\}\) such that $T_ny=\mu y$. By elliptic regularity, \(\phi_n\) is smooth up to the boundary, and hence \(\phi_n y\) belongs to the complexification of \(\mathbb X\). Let \(z(t):=T_n(t,0)y\). Then
\[
	z'(t)=\bigl(A(t)-\lambda_nD(t)\bigr)z(t).
\]
Since \(-\Delta\phi_n=\lambda_n\phi_n\), we have
\[
	\frac{\partial}{\partial t}\bigl(\phi_n(x)z(t)\bigr)
	=
	D(t)\Delta\bigl(\phi_n(x)z(t)\bigr)
	+
	A(t)\phi_n(x)z(t).
\]
Thus \(\phi_n(x)z(t)\) solves system \eqref{pthirdsys-linear-PDE}. By uniqueness of solutions, we obtain
\[
	\mathcal T\bigl[\phi_n y\bigr]
	=
	\phi_n T_ny
	=
	\mu\phi_n y .
\]
Hence \(\mu\in\sigma(\mathcal T)\). Therefore
\[
	\sigma(\mathcal T)\setminus\{0\}
	=
	\bigcup_{n\ge0}\sigma(T_n).
\]

We next show that modes are stable when $n$ is sufficiently large. Let $S_n(t,s)$ be the evolution operator associated with $v'=-\lambda_nD(t)v$. Since $D(t)$ is diagonal and $d_i(t)\ge d_*$, we have
\[
	\|S_n(t,s)\|\le e^{-\lambda_nd_*(t-s)},\qquad 0\le s\le t.
\]
For the solution $y(t)$ of system \eqref{pthirdsys-linear-decomposition-discrete} with $y(0)=y_0$, the variation-of-constants formula gives
\[
	y(t)=S_n(t,0)y_0+\int_0^t S_n(t,s)A(s)y(s)\,\mathrm ds.
\]
Therefore
\[
	\|y(t)\|
	\le
	e^{-\lambda_nd_*t}\|y_0\|
	+
	\int_0^t e^{-\lambda_nd_*(t-s)}A^*\|y(s)\|\,\mathrm ds.
\]
Setting $\psi(t):=e^{\lambda_nd_*t}\|y(t)\|$, we obtain
\[
	\psi(t)\le \|y_0\|+\int_0^t A^*\psi(s)\,\mathrm ds.
\]
By Gronwall's inequality,
\[
	\|y(t)\|\le \|y_0\|e^{(A^*-\lambda_nd_*)t}.
\]
Taking \(t=\omega\), we obtain
\begin{equation}\label{pthird-high-mode-estimate}
		\|T_n\|
		\le
		e^{\omega(A^*-\lambda_nd_*)}.
\end{equation}
Since \(\lambda_n\to+\infty\), the right-hand side is less than \(1\) for all sufficiently large \(n\). Therefore $\rho(T_n)\le \|T_n\|<1$ for all sufficiently large \(n\). This proves (ii).

It remains to prove (iii). By part~(i), $\mathcal T=D\mathcal P(\check u^*(0))$ has an eigenvalue of modulus greater than one. Since $\mathcal T$ is compact,
\cite[p.~107, Exercise~5]{Henry-book-1981} yields a nontrivial backward
orbit $\{\phi_{-k}\}_{k\ge0}$ such that
\[
\mathcal P(\phi_{-k})=\phi_{-k+1}\quad(k\ge1),
\qquad
\phi_{-k}\to\check u^*(0)\quad\text{as }k\to\infty.
\]

If the kinetic system~\eqref{pthirdsys-nonlinear-ODE} is
nonautonomous, then after reindexing, we may assume that $\phi_0\ne\check u^*(0)$. Since $u(k\omega,\phi_{-k})=\mathcal P^k(\phi_{-k})=\phi_0$ and
$\check u^*(k\omega)=\check u^*(0)$, we have $\|\phi_{-k}-\check u^*(0)\|_{\mathbb X}\to0$ and 
\[
	\|u(k\omega,\phi_{-k})-\check u^*(k\omega)\|_{\mathbb X}
	=\|\phi_0-\check u^*(0)\|_{\mathbb X}>0.
\]
Hence $\check u^*$ is Lyapunov unstable in $\mathbb X$.

If the kinetic system~\eqref{pthirdsys-nonlinear-ODE} is autonomous, then every point of $\check\Gamma$ is fixed by $\mathcal P$ since $
\mathcal P(\check u^*(t))
=\check u^*(t+\omega)
=\check u^*(t)$ for all $t\in\mathbb R$. The backward orbit cannot lie entirely in $\check\Gamma$. Otherwise, $\phi_{-k+1}=\mathcal P(\phi_{-k})=\phi_{-k}$ for all $k\ge1$, so it would be constant and hence identically equal to $\check u^*(0)$, contradicting its nontriviality.
After reindexing, we may therefore assume that
$\phi_0\notin\check\Gamma$. Since $\check\Gamma$ is compact, we have $\operatorname{dist}_{\mathbb X}(\phi_{-k},\check\Gamma)\to0$ and 
\[
	\operatorname{dist}_{\mathbb X}
	\bigl(u(k\omega,\phi_{-k}),\check\Gamma\bigr)
	=\operatorname{dist}_{\mathbb X}(\phi_0,\check\Gamma)>0.
\]
Thus $\check u^*$ is orbitally unstable in $\mathbb X$. This proves (iii).
\end{proof}

Combining the instability part of the preceding lemma with the kinetic stability assumption gives the following sufficient condition for Turing instability.

\begin{proposition}\label{pthirdprop-turing-instability}
Let (A1) and (A2) hold. If there exists \(n_0\ge1\) such that $\rho(T_{n_0})>1$, then a Turing instability occurs. Moreover, only finitely many spatial modes can satisfy the spectral instability condition \(\rho(T_n)>1\).
\end{proposition}

\begin{proof}
By assumption~(A2), \(u^*(t)\) is stable for the kinetic system. By Lemma~\ref{pthirdlem-instability-reduction}~(iii), \(\check{u}^*(t)\) is unstable for system~\eqref{pthirdsys-nonlinear-PDE} in \(\mathbb X\). Hence, by the definition of Turing instability, a Turing instability occurs. The final assertion follows from Lemma~\ref{pthirdlem-instability-reduction}(ii).
\end{proof}

Note that the Neumann eigenvalues depend on the domain \(\Omega\). To separate the analysis of the finite-dimensional modal equations from the geometry of a fixed domain, we introduce a continuous parameter \(\alpha\ge0\). For \(\alpha\ge0\), let \(T_\alpha(t,s)\) denote the evolution operator associated with
\begin{equation}\label{pthirdsys-linear-decomposition-continuous}
	\frac{\mathrm d y}{\mathrm d t}
	=
	\bigl(A(t)-\alpha D(t)\bigr)y,
\end{equation}
and set $T_\alpha:=T_\alpha(\omega,0)$. Then \(T_0\) is the monodromy matrix of the kinetic variational system.  Thus, on a fixed domain \(\Omega\), the discrete modal family is obtained by sampling the continuous family \(T_\alpha\) at the Neumann eigenvalues \(\alpha=\lambda_n\). In particular, \(\alpha=0\) corresponds to the spatially homogeneous mode, whereas a value \(\alpha>0\) corresponds to an actual spatially heterogeneous mode on \(\Omega\) precisely when \(\alpha\) is a nonzero Neumann eigenvalue. An arbitrary value \(\alpha>0\) can still be realized as a nonzero Neumann eigenvalue after a suitable uniform rescaling of the domain.

\begin{lemma}\label{pthirdlem-domain-scaling}
Assume that (A1) and (A2) hold, and  that \(\rho(T_\alpha)>1\) for some \(\alpha>0\).
Then the following statements are valid:
\begin{enumerate}
		\item[(i)]
		If \(\alpha=\lambda_n\) for some \(n\ge1\), then a Turing instability
		occurs for system~\eqref{pthirdsys-nonlinear-PDE} on the original
		domain \(\Omega\).
		
		\item[(ii)]
		Let \(\lambda_n>0\) be any nonzero Neumann eigenvalue of \(\Omega\). Set $\beta:=\sqrt{\frac{\lambda_n}{\alpha}}$ and $\widetilde\Omega:=\beta\Omega$. Then \(\alpha\) is a nonzero Neumann eigenvalue of \(\widetilde\Omega\), and a Turing instability occurs on \(\widetilde\Omega\).
\end{enumerate}
\end{lemma}

\begin{proof}
For part (i), the identity $T_\alpha=T_{\lambda_n}=T_n$ implies that \(\rho(T_n)>1\). The conclusion therefore follows directly from Proposition~\ref{pthirdprop-turing-instability}.
	
For part (ii), let \(\phi_n\) be a Neumann eigenfunction on \(\Omega\) associated with \(\lambda_n\), so that
\[
	-\Delta\phi_n=\lambda_n\phi_n,
	\qquad
	\frac{\partial\phi_n}{\partial\nu}=0.
\]
For \(\widetilde x\in\widetilde\Omega\), define $\widetilde\phi_n(\widetilde x):=\phi_n\left(\frac{\widetilde x}{\beta}\right)$. Writing \(x=\frac {\widetilde x}{\beta}\), we obtain
\[
	-\Delta_{\widetilde x}\widetilde\phi_n(\widetilde x)
	=
	\frac{\lambda_n}{\beta^2}
	\widetilde\phi_n(\widetilde x).
\]
Moreover, uniform dilation preserves the outward unit normal direction, and hence
\[
	\frac{\partial\widetilde\phi_n}{\partial\widetilde\nu}
	=
	\frac1\beta
	\frac{\partial\phi_n}{\partial\nu}
	=
	0
	\qquad\text{on }\partial\widetilde\Omega.
\]
Thus the corresponding Neumann eigenvalue on \(\widetilde\Omega\) is $\widetilde\lambda_n=\frac{\lambda_n}{\beta^2}=\alpha$. The assumptions (A1), (A2) and the kinetic periodic solution \(u^*\) are unaffected by replacing \(\Omega\) with \(\widetilde\Omega\). Therefore part {\rm(i)}, applied on \(\widetilde\Omega\), gives the desired Turing instability.
\end{proof}

Here we should point out that a similar continuous-parameter and domain-scaling argument was used in
\cite[Section~1]{Dancer-DIE-1994}.

\section{Nonoccurrence of Turing instability}\label{pthirdsec-Nonoccurrence-of-Turing-Instability}

In this section, we establish the estimate $\rho(T_n)<1$ for all
$n\ge1$ in three cases: equal diffusion rates, sufficiently large
diffusion rates on the fixed domain $\Omega$, and cooperative
systems. Under the additional hypotheses in
Theorem~\ref{pthirdthe-nonlinear-nonoccurrence}, these estimates imply
the nonoccurrence of Turing instability.

We shall repeatedly use the consequence $\rho(T_0)\le1$ of
assumption~(A2). Indeed, if $\rho(T_0)>1$, then the linearized
instability principle implies that $u^*$ is Lyapunov unstable
in the nonautonomous case and orbitally unstable in the
autonomous case, contradicting (A2).

\begin{lemma}\label{pthirdlem-diffusion-rate-equal}
Let (A1) and (A2) hold. If all the diffusion rates are equal, that is, $d_1(t)=d_2(t)=\dots =d_m(t)=:d(t)$ for all $t\in[0,\omega]$, then $\rho(T_n)<1$ for all $n\ge1$.
\end{lemma}

\begin{proof}
Since \(D(t)=d(t)I\), where $I$ is the identity matrix, system
\eqref{pthirdsys-linear-decomposition-discrete} becomes
\begin{equation}\label{pthirdsys-equal-rates}
	\frac{\mathrm d y}{\mathrm d t}
	=
	\big(A(t)-\lambda_n d(t)I\big)y.
\end{equation}
For any solution \(z(t)\) of the kinetic linearized system
$z'(t)=A(t)z(t)$, the function
\[
	y(t)
	=
	\exp\left(-\lambda_n \int_0^t d(s)\,\mathrm ds\right)z(t)
\]
solves system \eqref{pthirdsys-equal-rates}. Therefore the monodromy matrices satisfy
\[
	T_n
	=
	\exp\left(-\lambda_n\int_0^{\omega}d(s)\,\mathrm ds\right)T_0.
\]
It follows that
\[
	\rho(T_n)
	=
	\exp\left(-\lambda_n\int_0^{\omega}d(s)\,\mathrm ds\right)\rho(T_0).
\]
By assumption (A2), we have \(\rho(T_0)\le 1\). Hence
$\rho(T_n)<1$ for all $n>0$. 
\end{proof}

For $m=1$, there is only one diffusion coefficient, so
Lemma~\ref{pthirdlem-diffusion-rate-equal} applies directly to scalar equations. The next estimate uses the positive lower bound $\lambda_1>0$ on the fixed domain and does not require assumption~(A2).

\begin{lemma}\label{pthirdlem-diffusion-rate-large}
	Let (A1) hold. Then there exists \(\bar d>0\) such that, if $d_*\ge \bar d$, then $\rho(T_n)<1$ for all $n\ge1$.
\end{lemma}
\begin{proof}
Recall that $A^*:=\max_{t\in[0,\omega]}\|A(t)\|$. As in the proof of Lemma~\ref{pthirdlem-instability-reduction}, the evolution operator of the \(n\)-th mode satisfies
\[
	\|T_n(t,0)\|
	\le
	e^{(A^*-\lambda_n d_*)t},
	\qquad n\ge1.
\]
Taking \(t=\omega\), we obtain $\|T_n\|\le e^{\omega(A^*-\lambda_n d_*)}$ for all $n\ge1$. Since the domain is fixed, \(\lambda_n\ge\lambda_1>0\) for all \(n\ge1\). Choose \(\bar d>A^*/\lambda_1\). If \(d_*\ge\bar d\), then
\[
	A^*-\lambda_n d_*
	\le
	A^*-\lambda_1\bar d
	<0,
	\qquad n\ge1.
\]
Therefore
\[
	\|T_n\|
	\le
	e^{\omega(A^*-\lambda_n d_*)}
	<1,
	\qquad n\ge1.
\]
Hence
\[
	\rho(T_n)\le \|T_n\|<1,
	\qquad n\ge1.
\]
This completes the proof.
\end{proof}

We next turn from the size of the diffusion rates to the order structure of the
linearized reaction dynamics.

\begin{lemma}\label{pthirdlem-cooperative}
Assume that (A1) and (A2) hold, and that \(A(t)\) is cooperative for
all \(t\in[0,\omega]\). Then $\rho(T_n)<1$  for all $n\ge1$.
\end{lemma}

\begin{proof}
	
Let \(\mathbb R^m_+\) be the standard positive cone in \(\mathbb R^m\), and write \(x\ge y\) if \(x-y\in\mathbb R^m_+\). Since \(A(t)\) is cooperative, the matrices \(A(t)\) and \(A(t)-\lambda_nD(t)\) have nonnegative off-diagonal entries. Hence the corresponding evolution operators \(T_0(t,s)\) and \(T_n(t,s)\) are positive with respect to \(\mathbb R^m_+\) for \(t\ge s\) \cite[Proposition~2.1, p.~34]{Smith-book-1995}.

Fix \(n\ge1\) and \(x\in\mathbb R^m_+\), and set $u_n(t):=T_n(t,0)x$. Then \(u_n(t)\in\mathbb R^m_+\) for all \(t\ge0\). Let $v_n(t):=e^{\lambda_nd_*t}u_n(t)$. Since \(u_n\) solves \eqref{pthirdsys-linear-decomposition-discrete}, we have
\[
v_n'
=
\big(A(t)-\lambda_n(D(t)-d_*I)\big)v_n.
\]
Because \(D(t)-d_*I\) is a nonnegative diagonal matrix and
\(v_n(t)\in\mathbb R^m_+\), it follows componentwise that $v_n'(t)\le A(t)v_n(t)$. By the comparison theorem for cooperative linear systems, it follows that 
\[
0\le v_n(t)\le T_0(t,0)x,
\qquad t\ge0.
\]
Therefore
\[
0\le T_n(t,0)x
\le
e^{-\lambda_nd_*t}T_0(t,0)x,
\qquad x\in\mathbb R^m_+.
\]
Taking \(t=\omega\), we obtain the operator inequality
\[
0\le T_n
\le
e^{-\omega\lambda_nd_*}T_0
\]
with respect to the order induced by \(\mathbb R^m_+\). By the monotonicity
of the spectral radius for positive linear operators \cite[Lemma~2.2]{Lemmens-Nussbaum-PAMS-2013}, we have 
\[
\rho(T_n)
\le
e^{-\omega\lambda_nd_*}\rho(T_0).
\]
Note that Assumption~(A2) implies that \(\rho(T_0)\le1\). Hence
\[
\rho(T_n)
\le
e^{-\omega\lambda_nd_*}
<1,
\qquad n\ge1.
\]
This completes the proof.
\end{proof}

\begin{remark}\label{pthirdrem-general-cone}
Under (A1) and (A2), the conclusion of
Lemma~\ref{pthirdlem-cooperative} remains valid if there exists a proper cone \(K\subset\mathbb R^m\) and a constant \(\eta>0\) such that
\[
	T_0(t,s)K\subset K,\qquad
	T_n(t,s)K\subset K,\qquad
	(D(t)-\eta I)K\subset K,
	\qquad t\ge s.
\]
Then
\[
0\le_K T_n(\omega,0)
\le_K
e^{-\omega\lambda_n\eta}T_0(\omega,0),
\]
and therefore
\[
\rho(T_n)
\le
e^{-\omega\lambda_n\eta}\rho(T_0)
<1,
\qquad n\ge1.
\]
Thus, no nonzero Neumann mode can satisfy the condition \(\rho(T_n)>1\). See, for example, \cite[Section~3.5]{Smith-book-1995} for background on monotone dynamical systems and 
cone-induced orders.
\end{remark}

We now combine the preceding spectral estimates with the linearized stability
principle for the period map to obtain nonlinear nonoccurrence results.

\begin{theorem}\label{pthirdthe-nonlinear-nonoccurrence}
Let  (A1) and (A2) hold, and assume that all diffusion rates are equal, that $d_*\ge\bar d$ with $\bar d$ as in Lemma~\ref{pthirdlem-diffusion-rate-large}, or that $A(t)$ is cooperative for all $t\in[0,\omega]$. If one of the following conditions holds:
\begin{enumerate}
		\item[(i)] $\rho(T_0)<1$;
		\item[(ii)] system~\eqref{pthirdsys-nonlinear-PDE} is autonomous, $u^*(t)$ is nonconstant, $1$ is an algebraically simple Floquet multiplier of $T_0$, and all other Floquet
		multipliers of $T_0$ lie strictly inside the unit circle,
\end{enumerate}
then Turing instability does not occur.
\end{theorem}

\begin{proof}
By Lemmas~\ref{pthirdlem-diffusion-rate-equal}, \ref{pthirdlem-diffusion-rate-large}, and \ref{pthirdlem-cooperative}, we have $\rho(T_n)<1$ for all $n\ge1$. In case (i), Lemma~\ref{pthirdlem-instability-reduction}~(i)--(ii), together with $\rho(T_0)<1$, gives $\rho(\mathcal T)<1$. Since $\mathcal T=D\mathcal P(\check u^*(0))$, the principle of linearized stability implies that $\check u^*(0)$ is an asymptotically stable fixed point of $\mathcal P$ \cite[p.~100, Exercise~2]{Henry-book-1981}. By periodicity and continuous dependence over one period, $\check u^*(t)$ is asymptotically stable for system~\eqref{pthirdsys-nonlinear-PDE} in $\mathbb X$.

	In case (ii), $1\notin\sigma(T_n)$ for every $n\ge1$.
Fix $k\ge1$ and let $\xi$ belong to the complexification of
$\mathbb X$ with $(\mathcal T-I)^k\xi=0$.
The Fourier coefficient identity in the proof of
Lemma~\ref{pthirdlem-instability-reduction}(i) gives
\[
0
=
c_n\bigl((\mathcal T-I)^k\xi\bigr)
=
(T_n-I)^k c_n(\xi),
\qquad n\ge0.
\]
Since $T_n-I$ is invertible for $n\ge1$, we have
$c_n(\xi)=0$ for all $n\ge1$. By completeness of the Neumann
eigenfunctions, $\xi$ is spatially homogeneous.
On the spatially homogeneous subspace, $\mathcal T$ is identified
with $T_0$. Consequently,
\[
\ker\bigl((\mathcal T-I)^k\bigr)
=
\left\{
\check a:
a\in\ker\bigl((T_0-I)^k\bigr)
\right\},
\qquad k\ge1.
\]
Since $1$ is algebraically simple for $T_0$, it is also
algebraically simple for $\mathcal T$.
Moreover, $\partial_t\check u^*$ is a nonzero $\omega$-periodic
solution of the linearized equation. It then follows that 
\[
\mathcal T\,\partial_t\check u^*(0)
=
\partial_t\check u^*(\omega)
=
\partial_t\check u^*(0)\ne0,
\]
and hence, in the real phase space,
\[
\ker(\mathcal T-I)
=
\operatorname{span}\{\partial_t\check u^*(0)\}.
\]
By compactness of $\mathcal T$, the multiplier $1$ is isolated.
Lemma~\ref{pthirdlem-instability-reduction}(i)--(ii) and the
spectral assumptions further imply that, for some
$\epsilon\in(0,1)$,
\[
\sigma(\mathcal T)\setminus\{1\}
\subset
\{z\in\mathbb C:|z|\le\epsilon\}.
\]
Since system~\eqref{pthirdsys-nonlinear-PDE} is autonomous,
\cite[Theorem~8.2.3]{Henry-book-1981} implies that
$\check\Gamma$ is orbitally asymptotically stable for
system~\eqref{pthirdsys-nonlinear-PDE} in $\mathbb X$.

Thus, $\check{u}^*(t)$ is stable for
system~\eqref{pthirdsys-nonlinear-PDE}  in both cases.
Since $u^*(t)$ is stable for the kinetic system by (A2), it follows that 
Turing instability does not occur.
\end{proof}


\section{Criteria for Turing instability}\label{pthirdsec-sufficient-conditions}
In this section, we seek verifiable sufficient conditions for Turing instability in periodic reaction--diffusion systems. To this end, we start with a class of periodic Riccati equations.

\subsection{A class of Riccati equations}
\label{pthirdsubsec-periodic-riccati}

Consider the following periodic Riccati equation:
\begin{equation}\label{pthirdsys-general-Riccati-alpha}
	W'
	=
	g(t)\bigl(1+W^2\bigr)
	-
	\bigl[h(t)+\alpha\kappa(t)\bigr]W,
	\qquad \alpha\in\mathbb R.
\end{equation}
Throughout this subsection, assume that 
$g,h,\kappa$ are continuous and $\omega$-periodic on $\mathbb R$, and that 
$g(t)>0$ and $\kappa(t)>0$ for all $t\in\mathbb R$. Set
\[
F_\alpha(t,W)
:=
g(t)(1+W^2)-\bigl[h(t)+\alpha\kappa(t)\bigr]W.
\]
For each fixed $\alpha\in\mathbb R$ and each $\zeta\ge0$, let
$W(t;\alpha,\zeta)$ denote the solution of
\eqref{pthirdsys-general-Riccati-alpha} satisfying
$W(0;\alpha,\zeta)=\zeta$, whenever it is defined. Define
\[
\mathcal J_\alpha
:=
\left\{
\zeta\ge0:
W(t;\alpha,\zeta)
\text{ exists on }[0,\omega]
\right\}.
\]
For $\zeta\in\mathcal J_\alpha$, define the period map
$\widetilde P_\alpha(\zeta):=W(\omega;\alpha,\zeta)$. We shall also use
$\mathcal I_\alpha:=\mathcal J_\alpha\cap(0,\infty)$. For
$\zeta\in\mathcal I_\alpha$, set
\begin{equation}\label{pthird-fixpoint-to-zeros}
	G_\alpha(\zeta):=\widetilde P_\alpha(\zeta)-\zeta.
\end{equation}
Then the existence of a positive $\omega$-periodic solution of the
Riccati equation \eqref{pthirdsys-general-Riccati-alpha} is equivalent
to the existence of a positive zero of $G_\alpha$. Indeed, if
$\zeta\in\mathcal I_\alpha$ and $G_\alpha(\zeta)=0$, then
$W(\omega;\alpha,\zeta)=W(0;\alpha,\zeta)$. Since the right-hand side
of \eqref{pthirdsys-general-Riccati-alpha} is $\omega$-periodic in $t$,
uniqueness of solutions implies
$W(t+\omega;\alpha,\zeta)=W(t;\alpha,\zeta)$ for all $t\in\mathbb R$.
Moreover, since $\zeta>0$ and $F_\alpha(t,0)=g(t)>0$, this solution
is positive for all $t\in\mathbb R$. Conversely, every positive
$\omega$-periodic solution $W(t)$ gives a zero of $G_\alpha$ by taking
$\zeta=W(0)$.

We now classify the positive $\omega$-periodic solutions of
\eqref{pthirdsys-general-Riccati-alpha} by means of the Riccati period
map. Set
\[
\Theta(\alpha)
:=
\inf_{\zeta\in\mathcal I_\alpha}G_\alpha(\zeta),
\]
with the convention that $\Theta(\alpha)=+\infty$ if
$\mathcal I_\alpha=\emptyset$.

\begin{lemma}\label{pthirdthm-fixed-alpha-classification}
Let  $\alpha\in\mathbb R$ be given. Then the following statements are valid:
\begin{enumerate}
		\item[(i)]
		If $\Theta(\alpha)>0$, then equation
		\eqref{pthirdsys-general-Riccati-alpha} has no positive
		periodic solution.
		
		\item[(ii)]
		If $\Theta(\alpha)=0$, then equation
		\eqref{pthirdsys-general-Riccati-alpha} has a unique positive
		periodic solution $W^*$.
		
		\item[(iii)]
		If $\Theta(\alpha)<0$, then equation
		\eqref{pthirdsys-general-Riccati-alpha} has exactly two
		distinct positive periodic solutions $W_{\alpha,-}$ and
		$W_{\alpha,+}$ with $W_{\alpha,-}(t)<W_{\alpha,+}(t)$ for all
		$t\in\mathbb R$.
\end{enumerate}
\end{lemma}

\begin{proof}
Fix $\alpha\in\mathbb R$. Since $g(t)>0$, we have $F_\alpha(t,0)=g(t)>0$. Hence a solution starting from a positive initial value cannot cross the line $W=0$ from above. Therefore every solution with $\zeta>0$ remains positive as long as it exists.

We first describe the interval $\mathcal I_\alpha$. By definition, $\mathcal I_\alpha$ consists of those positive initial values for which the corresponding solution is defined on the whole interval $[0,\omega]$. By continuous dependence on initial data, $\mathcal I_\alpha$ is open in $(0,\infty)$. Moreover, if $0<\zeta_1<\zeta_2$ and both corresponding solutions exist on a
common interval, then $W(t;\alpha,\zeta_1)<W(t;\alpha,\zeta_2)$ on that interval. Consequently, if $\zeta_2\in\mathcal I_\alpha$, then $\zeta_1\in\mathcal I_\alpha$. Thus either $\mathcal I_\alpha=\emptyset$, or $\mathcal I_\alpha=(0,R_\alpha)$ for some $R_\alpha\in(0,\infty]$. We now show that, whenever $\mathcal I_\alpha\neq\emptyset$, we actually has $R_\alpha<\infty$. Let
\[
	g_0:=\min_{t\in[0,\omega]}g(t)>0,
	\qquad
	B^*:=\max_{t\in[0,\omega]}|h(t)+\alpha\kappa(t)|.
\]
Choose $M_0>0$ sufficiently large so that
\[
	g_0W^2-B^*W+g_0
	\ge
	\frac{g_0}{2}W^2,
	\qquad \forall~ W\ge M_0.
\]
Then, for $W\ge M_0$, we have
\[
	F_\alpha(t,W)
	\ge
	\frac{g_0}{2}W^2.
\]
Now choose $\zeta>\max\left\{M_0,\frac{2}{g_0\omega}\right\}$ and write $W(t):=W(t;\alpha,\zeta)$. Since $W(0)=\zeta\ge M_0$ and $W'(t)\ge\frac{g_0}{2}W(t)^2>0$ whenever $W(t)\ge M_0$, the solution satisfies $W(t)\ge\zeta\ge M_0$ throughout its maximal interval of existence. Consequently,
\[
	\left(\frac{1}{W(t)}\right)'
	=
	-\frac{W'(t)}{W(t)^2}
	\le
	-\frac{g_0}{2}.
\]
Integrating from $0$ to $t$, we obtain $\frac{1}{W(t)}\le\frac{1}{\zeta}-\frac{g_0}{2}t$. Since $\zeta>\frac{2}{g_0\omega}$, the right-hand side becomes nonpositive before $t=\omega$. This is impossible as long as the solution remains finite and positive. Hence the solution starting from such a $\zeta$ cannot be defined on the whole interval $[0,\omega]$. Therefore $\zeta\notin\mathcal I_\alpha$.
Since this holds for all sufficiently large $\zeta$, we conclude that $R_\alpha<\infty$.
	
Assume now that $\mathcal I_\alpha\neq\emptyset$. Since a positive solution with initial value in $\mathcal I_\alpha$ exists on $[0,\omega]$, the solution with initial value $W(0)=0$ also exists on $[0,\omega]$ and is dominated by it. Hence $0\in\mathcal J_\alpha$, and therefore $\widetilde P_\alpha(0)=W(\omega;\alpha,0)$. Since $F_\alpha(t,0)=g(t)>0$, this solution satisfies $W(t;\alpha,0)>0$ for $t\in(0,\omega]$. In particular,
\begin{equation}\label{pthirdineq-near-zero}
		\lim_{\zeta\to0^+}G_\alpha(\zeta)
		=
		\widetilde P_\alpha(0)>0.
\end{equation}
Next we prove the behavior at the right endpoint $R_\alpha$. Since $\widetilde P_\alpha$ is increasing, the limit $L:=\lim_{\zeta\to R_\alpha^-}\widetilde P_\alpha(\zeta)$ exists in $(0,\infty]$. We claim that $L=+\infty$. Suppose, by
contradiction, that $L<+\infty$. Let $H:=\max\{M_0,L\}$. Then, for every $\zeta\in(0,R_\alpha)$ and every $t\in[0,\omega]$, we must have $W(t;\alpha,\zeta)\le H$. Indeed, if $W(t_0;\alpha,\zeta)>H$ for some $t_0\in[0,\omega]$, then $W(t;\alpha,\zeta)\ge M_0$ for $t\ge t_0$ as long as the solution exists, and hence $W'(t)>0$. It follows that $W(\omega;\alpha,\zeta)>H\ge L$, which contradicts the definition of $L$. Thus the family $W(t;\alpha,\zeta)$ is uniformly bounded on $[0,\omega]$ as $\zeta\to R_\alpha^-$. By continuation and continuous dependence, the solution with initial value $W(0)=R_\alpha$ also exists on $[0,\omega]$. This contradicts the definition of $R_\alpha$ as the right endpoint of $\mathcal I_\alpha$. Therefore, $\lim_{\zeta\to R_\alpha^-}\widetilde P_\alpha(\zeta)=+\infty$, and hence
\begin{equation}\label{pthirdineq-near-R_alpha}
		\lim_{\zeta\to R_\alpha^-}G_\alpha(\zeta)=+\infty.
\end{equation}
	
We now prove the strict convexity of $G_\alpha$. By the differentiable dependence of solutions on initial data, for every $\zeta\in\mathcal I_\alpha$,
\[
	\partial_\zeta W(t;\alpha,\zeta)
	=
	\exp\left\{
	\int_0^t
	\partial_WF_\alpha
	\bigl(s,W(s;\alpha,\zeta)\bigr)\,\mathrm ds
	\right\}
	>0.
\]
In particular, $\widetilde P_\alpha'(\zeta)=\partial_\zeta W(\omega;\alpha,\zeta)>0$. Differentiating this expression with respect to $\zeta$,
we obtain
\[
	\widetilde P_\alpha''(\zeta)
	=
	\widetilde P_\alpha'(\zeta)
	\int_0^\omega
	\partial_{WW}F_\alpha
	\bigl(t,W(t;\alpha,\zeta)\bigr)
	\partial_\zeta W(t;\alpha,\zeta)\,\mathrm dt.
\]
Since $\partial_{WW}F_\alpha(t,W)=2g(t)>0$ and $\partial_\zeta W(t;\alpha,\zeta)>0$, it follows that $\widetilde P_\alpha''(\zeta)>0$ for all $\zeta\in\mathcal I_\alpha$. Consequently,
\[
	G_\alpha''(\zeta)
	=
	\widetilde P_\alpha''(\zeta)>0,
	\qquad
	\forall \zeta\in\mathcal I_\alpha,
\]
and hence $G_\alpha$ is strictly convex on $\mathcal I_\alpha$.
	
By the correspondence between positive $\omega$-periodic solutions and zeros of $G_\alpha$, it remains only to classify the zeros of $G_\alpha$. If $\mathcal I_\alpha=\emptyset$, then by definition there is no positive $\omega$-periodic solution and $\Theta(\alpha)=+\infty>0$. Hence the classification is immediate in this case. Assume therefore that $\mathcal I_\alpha=(0,R_\alpha)$. We have established inequalities \eqref{pthirdineq-near-zero} and \eqref{pthirdineq-near-R_alpha}, and shown that $G_\alpha$ is strictly convex. If $\Theta(\alpha)>0$, then $G_\alpha(\zeta)>0$ for all $\zeta\in\mathcal I_\alpha$, so $G_\alpha$ has no zero. Therefore equation \eqref{pthirdsys-general-Riccati-alpha} has no positive $\omega$-periodic solution. If $\Theta(\alpha)=0$, then the endpoint behavior above implies that the infimum is attained at some interior point $\zeta_\ast\in\mathcal I_\alpha$. Thus $G_\alpha(\zeta_\ast)=0$. Since $G_\alpha$ is strictly convex, this zero is unique. Hence equation \eqref{pthirdsys-general-Riccati-alpha} has exactly one positive $\omega$-periodic solution $W^*$. Finally, suppose $\Theta(\alpha)<0$. Then there exists $\zeta_0\in\mathcal I_\alpha$ such that $G_\alpha(\zeta_0)<0$. Together with inequalities \eqref{pthirdineq-near-zero} and \eqref{pthirdineq-near-R_alpha}, the intermediate value theorem gives at least two zeros of $G_\alpha$, one on each side of $\zeta_0$. Strict convexity implies that $G_\alpha$ cannot have more than two zeros. Therefore equation
\eqref{pthirdsys-general-Riccati-alpha} has exactly two positive $\omega$-periodic solutions. Denote the corresponding initial values by $0<\zeta_{\alpha,-}<\zeta_{\alpha,+}<R_\alpha$, and define $W_{\alpha,\pm}(t):=W(t;\alpha,\zeta_{\alpha,\pm})$. By the order-preserving property of scalar ordinary differential equations, $W_{\alpha,-}(t)<W_{\alpha,+}(t)$ for all $t\in[0,\omega]$. Since both solutions are $\omega$-periodic, the same inequality holds for all $t\in\mathbb R$.
\end{proof}

Let
\[
\mathcal A
:=
\left\{
\alpha\in\mathbb R:
\text{equation} \eqref{pthirdsys-general-Riccati-alpha}
\text{ admits a positive $\omega$-periodic solution}
\right\},
\]
and define $\alpha_*=\inf\mathcal A$.

\begin{lemma}\label{pthirdlem-alpha-set-structure}
 $\mathcal A=[\alpha_*,\infty)$, and $\alpha_*$
satisfies
\[
	\frac{
		\displaystyle\int_0^\omega
		\left[2g(t)-h(t)\right]\,\mathrm dt
	}{
		\displaystyle\int_0^\omega\kappa(t)\,\mathrm dt
	}
	\le
	\alpha_*
	\le
	\max_{t\in[0,\omega]}
	\frac{2g(t)-h(t)}{\kappa(t)}.
\]
\end{lemma}

\begin{proof}
We first prove that $G_\alpha$ is strictly decreasing with respect to $\alpha$. Let $\alpha_1<\alpha_2$, and let $\zeta\in\mathcal I_{\alpha_1}$. Denote by $W_i(t):=W(t;\alpha_i,\zeta)$, $i=1,2$, the corresponding solutions, whenever they are defined. Since
\[
	F_{\alpha_2}(t,W)
	=
	F_{\alpha_1}(t,W)
	-
	(\alpha_2-\alpha_1)\kappa(t)W,
\]
and since $\kappa(t)>0$, we have $F_{\alpha_2}(t,W)<F_{\alpha_1}(t,W)$ for all $W>0$. Moreover, positive solutions cannot cross the lower boundary $W=0$ from above because $F_\alpha(t,0)=g(t)>0$. Therefore, the solution $W_2(t)$ exists on $[0,\omega]$ and satisfies $0<W_2(t)<W_1(t)$ on $0<t\le\omega$. Consequently, $\mathcal I_{\alpha_1}\subset\mathcal I_{\alpha_2}$ and $P_{\alpha_2}(\zeta)<P_{\alpha_1}(\zeta)$ for all $\zeta\in\mathcal I_{\alpha_1}$. Equivalently,
\begin{equation}\label{pthird-alpha-monotonicity-G}
		G_{\alpha_2}(\zeta)<G_{\alpha_1}(\zeta),
		\qquad
		\forall \zeta\in\mathcal I_{\alpha_1}.
\end{equation}
	
We next show that $\mathcal A$ is nonempty. Set
\[
	\overline\alpha
	:=
	\max_{t\in[0,\omega]}
	\frac{2g(t)-h(t)}{\kappa(t)}.
\]
This number is finite since $\kappa(t)>0$ on $[0,\omega]$. If $\alpha\ge\overline\alpha$, then
\[
	h(t)+\alpha\kappa(t)\ge2g(t),
	\qquad t\in[0,\omega].
\]
Hence, for the constant function $\overline W(t)\equiv1$,
\[
	\overline W'(t)=0
	\ge
	g(t)(1+\overline W^2(t))
	-
	[h(t)+\alpha\kappa(t)]\overline W(t).
\]
Thus the interval $0\le W\le1$ is positively invariant for the equation \eqref{pthirdsys-general-Riccati-alpha}. The Poincar\'e map $\widetilde P_\alpha$ is therefore well-defined on $[0,1]$ and satisfies $\widetilde P_\alpha([0,1])\subset[0,1]$. Furthermore, we have $\widetilde P_\alpha(0)>0$ and $\widetilde P_\alpha(1)\le1$. By the intermediate value theorem, there exists $\zeta_\alpha\in(0,1]$ such that $\widetilde P_\alpha(\zeta_\alpha)=\zeta_\alpha$. Therefore equation \eqref{pthirdsys-general-Riccati-alpha} admits a positive $\omega$-periodic solution. Hence $\mathcal A\neq\emptyset$ and
\[
	\alpha_*
	\le
	\overline\alpha
	=
	\max_{t\in[0,\omega]}
	\frac{2g(t)-h(t)}{\kappa(t)}.
\]
We now prove that $\mathcal A$ is bounded from below. Let $\alpha\in\mathcal A$, and let $W_\alpha(t)$ be a positive $\omega$-periodic solution of equation
\eqref{pthirdsys-general-Riccati-alpha}. Dividing the Riccati equation by $W_\alpha(t)>0$ gives
\[
	\frac{W_\alpha'(t)}{W_\alpha(t)}
	=
	g(t)\left[
	W_\alpha(t)+\frac1{W_\alpha(t)}
	\right]
	-
	h(t)
	-
	\alpha\kappa(t).
\]
Integrating over one period and using the periodicity of $W_\alpha$ gives
\[
	0
	=
	\int_0^\omega
	g(t)\left[
	W_\alpha(t)+\frac1{W_\alpha(t)}
	\right]\,\mathrm dt
	-
	\int_0^\omega[h(t)]\,\mathrm dt
	-
	\alpha\int_0^\omega\kappa(t)\,\mathrm dt.
\]
Since $W_\alpha(t)+\frac1{W_\alpha(t)}\ge2$, we obtain
\[
	\alpha
	\ge
	\frac{
		\displaystyle\int_0^\omega
		\left[2g(t)-h(t)\right]\,\mathrm dt
	}{
		\displaystyle\int_0^\omega\kappa(t)\,\mathrm dt
	}.
\]
Thus, $\mathcal A$ is bounded from below, and hence, $\alpha_*$ is finite.
	
We next show that $\mathcal A$ is upward closed, that is, if $\beta_1\in\mathcal A$ and $\beta_2>\beta_1$, then $\beta_2\in\mathcal A$. Let $W_{\beta_1}(t)$ be a positive $\omega$-periodic solution of \eqref{pthirdsys-general-Riccati-alpha} with $\alpha=\beta_1$, and set $\zeta_{\beta_1}:=W_{\beta_1}(0)$. Then $P_{\beta_1}(\zeta_{\beta_1})=\zeta_{\beta_1}$. By the strict monotonicity
\eqref{pthird-alpha-monotonicity-G}, we have $P_{\beta_2}(\zeta_{\beta_1})<P_{\beta_1}(\zeta_{\beta_1})=\zeta_{\beta_1}$.
Hence $G_{\beta_2}(\zeta_{\beta_1})<0$. By Lemma~\ref{pthirdthm-fixed-alpha-classification}, this implies that \eqref{pthirdsys-general-Riccati-alpha} admits positive $\omega$-periodic solutions for this $\alpha=\beta_2$. Therefore $\beta_2\in\mathcal A$. It follows that $(\alpha_*,\infty)\subset\mathcal A$.

It remains to prove that $\alpha_*\in\mathcal A$. Choose a decreasing sequence $\alpha_n\rightarrow\alpha_*$ with $\alpha_n>\alpha_*$. By the previous paragraph, $\alpha_n\in\mathcal A$ for every $n$. Let $W_n(t)$ be a positive $\omega$-periodic solution of equation \eqref{pthirdsys-general-Riccati-alpha} with $\alpha=\alpha_n$. Since $\alpha_n\to\alpha_*$, there exists $M>0$ such that $\|h(t)+\alpha_n\kappa(t)\|_\infty\le M$ for all $n$. Recall that $g_0:=\min_{t\in[0,\omega]}g(t)>0$. At every maximum or minimum point $t_n$ of the periodic function $W_n$, we have $W_n'(t_n)=0$, and hence
\[
	g(t_n)\left(
	W_n(t_n)+\frac1{W_n(t_n)}
	\right)
	=
	h(t_n)+\alpha_n\kappa(t_n).
\]
Consequently,
\[
	W_n(t_n)+\frac1{W_n(t_n)}
	\le
	\frac{M}{g_0}=:K.
\]
Since the left-hand side is at least $2$, we have $K\ge2$. Applying this estimate at a maximum point and at a minimum point of $W_n$, respectively, gives
\[
	K^{-1}\le W_n(t)\le K,
	\qquad t\in[0,\omega].
\]
The Riccati equation then yields a uniform bound on $W_n'$. Hence, by the Arzel\`a--Ascoli theorem, after passing to a subsequence, we may assume that $W_n\to W_\infty$ uniformly on $[0,\omega]$. Passing to the limit in the integral equation
\[
	W_n(t)
	=
	W_n(0)
	+
	\int_0^t
	\left[
	g(s)(1+W_n^2(s))
	-
	(h(s)+\alpha_n\kappa(s))W_n(s)
	\right]\,\mathrm ds,
\]
we obtain
\[
	W_\infty(t)
	=
	W_\infty(0)
	+
	\int_0^t
	\left[
	g(s)(1+W_\infty^2(s))
	-
	(h(s)+\alpha_*\kappa(s))W_\infty(s)
	\right]\,\mathrm ds.
\]
Thus, $W_\infty$ is a solution of equation \eqref{pthirdsys-general-Riccati-alpha} with $\alpha=\alpha_*$. Moreover,
\[
	K^{-1}\le W_\infty(t)\le K,
	\qquad t\in[0,\omega],
\]
so $W_\infty$ is positive. Since each $W_n$ is $\omega$-periodic, the uniform limit $W_\infty$ is also $\omega$-periodic. Hence $\alpha_*\in\mathcal A$. By the  upward closedness of $\mathcal A$ 
and $\alpha_*=\inf\mathcal A$, it then follows  that $\mathcal A=[\alpha_*,\infty)$.
\end{proof}

\begin{lemma}\label{pthirdthm-alpha-critical-classification}
For every $\alpha\in\mathbb R$, $\Theta(\alpha)$ and $\alpha_*-\alpha$ have the same sign.
\end{lemma}

\begin{proof}
By Lemma~\ref{pthirdlem-alpha-set-structure}, we have $\mathcal A=[\alpha_*,\infty)$. If $\alpha<\alpha_*$, then $\alpha\notin\mathcal A$, and
Lemma~\ref{pthirdthm-fixed-alpha-classification} gives $\Theta(\alpha)>0$.
	
Since $\alpha_*\in\mathcal A$, we have $\Theta(\alpha_*)\le0$. Suppose that $\Theta(\alpha_*)<0$. Then there exists $\zeta_0\in\mathcal I_{\alpha_*}$ such that $G_{\alpha_*}(\zeta_0)<0$. By continuous dependence of solutions on parameters, for all $\alpha$ sufficiently close to $\alpha_*$, the solution with initial value $\zeta_0$ exists on $[0,\omega]$ and $G_\alpha(\zeta_0)<0$. In particular, this holds for some $\alpha<\alpha_*$. Lemma~\ref{pthirdthm-fixed-alpha-classification} then implies that $\alpha\in\mathcal A$, a contradiction. Therefore $\Theta(\alpha_*)=0$.

Finally, let $\alpha>\alpha_*$, and let $\zeta_*$ be a positive zero of $G_{\alpha_*}$. By \eqref{pthird-alpha-monotonicity-G},it follows that $G_\alpha(\zeta_*)<G_{\alpha_*}(\zeta_*)=0$, and hence, $\Theta(\alpha)<0$. 
\end{proof}

Now we are able to classify the global  dynamics of the Riccati equation \eqref{pthirdsys-general-Riccati-alpha}
 in terms of parameter $\alpha$.

\begin{theorem}\label{pthirdprop-riccati-global-dynamics}
The following statements are valid for equation \eqref{pthirdsys-general-Riccati-alpha}:
\begin{enumerate}
	\item[(i)]
		If $\alpha<\alpha_*$, then equation \eqref{pthirdsys-general-Riccati-alpha} has no positive
		$\omega$-periodic solution, and every positive solution
		blows up to $+\infty$ in finite forward time.
		
	\item[(ii)]
		If $\alpha=\alpha_*$, then equation
		\eqref{pthirdsys-general-Riccati-alpha} has a unique positive
		$\omega$-periodic solution $W^*$ and $\widetilde P_{\alpha_*}'(W^*(0))=1$.
		Moreover, for any $0<\zeta<W^*(0)$, the solution $W(t;\alpha_*,\zeta)$ of \eqref{pthirdsys-general-Riccati-alpha} with
	    $W(0)=\zeta$ exists globally on $[0,\infty]$ and $\lim_{t\to\infty}
		|W(t;\alpha_*,\zeta)-W^*(t)|=0$; for any $\zeta>W^*(0)$, the solution $W(t;\alpha_*,\zeta)$ blows up to $+\infty$ in finite time.
		
	\item[(iii)]
		If $\alpha>\alpha_*$, then equation
		\eqref{pthirdsys-general-Riccati-alpha} has exactly two
		positive $\omega$-periodic solutions $W_{\alpha,-}$ and
		$W_{\alpha,+}$ with $W_{\alpha,-}(t)<W_{\alpha,+}(t)$
		for all $t\in\mathbb R$, and $0<\widetilde P_\alpha'(W_{\alpha,-}(0))<1
		<\widetilde P_\alpha'(W_{\alpha,+}(0))$. Further, for any $0<\zeta<W_{\alpha,+}(0)$, the solution $W(t;\alpha_*,\zeta)$ of \eqref{pthirdsys-general-Riccati-alpha} with
		 $W(0)=\zeta$ exists globally on $[0,\infty]$ and $
		\lim_{t\to\infty}|W(t;\alpha,\zeta)-W_{\alpha,-}(t)|=0$;
		for any $\zeta>W_{\alpha,+}(0)$, the solution $W(t;\alpha,\zeta)$ blows up to $+\infty$ in finite time.
\end{enumerate}
\end{theorem}

\begin{proof}
By Lemmas~\ref{pthirdthm-fixed-alpha-classification} and \ref{pthirdthm-alpha-critical-classification}, it remains to prove the assertions about stability and forward dynamics. If $\zeta\notin\mathcal I_\alpha$, the positive solution blows up at or before $t=\omega$. Thus, whenever period-map iterates are used below, we may assume that $\zeta\in\mathcal I_\alpha$.
	
We first recall the multiplier of the Riccati period map at a positive periodic solution. Let $W_p(t)=W(t;\alpha,\zeta_p)$ be a positive $\omega$-periodic solution. The variational equation gives
\[
	\widetilde P_\alpha'(\zeta_p)
	=
	\exp\left\{
	\int_0^\omega
	\partial_WF_\alpha(t,W_p(t))\,\mathrm dt
	\right\}
	>0.
\]
Moreover, $G_\alpha$ is strictly convex on $\mathcal I_\alpha$, and positive $\omega$-periodic solutions are precisely the positive zeros of $G_\alpha$.

If $\alpha<\alpha_*$, then $\Theta(\alpha)>0$, and $G_\alpha(\zeta)>0$ for all $\zeta\in\mathcal I_\alpha$. Hence $\widetilde P_\alpha(\zeta)>\zeta$ on $\mathcal I_\alpha$. Let $\zeta_n:=\widetilde P_\alpha^n(\zeta)$, as long as the iterates are defined and belong to $\mathcal I_\alpha$. Then $\zeta_n$ is strictly increasing. If all iterates stayed in $\mathcal I_\alpha$, then either $\zeta_n$ would converge to an interior fixed point of $\widetilde P_\alpha$, which is impossible, or $\zeta_n\to R_\alpha$. The latter is also impossible since $\lim_{\zeta\to R_\alpha^-}\widetilde P_\alpha(\zeta)=+\infty$. So for $\zeta$ sufficiently close to $R_\alpha$, we have $\widetilde P_\alpha(\zeta)>R_\alpha$. Therefore some iterate leaves $\mathcal I_\alpha$, which means that the corresponding positive solution blows up during the following period. This proves (i).
	
If $\alpha=\alpha_*$, then $\Theta(\alpha)=0$. Let $\zeta_*=W^*(0)$ be the unique zero of $G_\alpha$. Since $G_\alpha$ is strictly convex and the endpoint
limits are positive, $\zeta_*$ is the unique interior minimizer of $G_\alpha$. Hence $G_\alpha'(\zeta_*)=0$, and therefore $P_\alpha'(\zeta_*)=1$. Moreover,
\[
	G_\alpha(\zeta)>0,
	\qquad
	\zeta\in\mathcal I_\alpha,\quad\zeta\neq\zeta_*.
\]
Thus $\widetilde P_\alpha(\zeta)>\zeta$ for $\zeta\neq\zeta_*$. Since $\widetilde P_\alpha$ is increasing, if $0<\zeta<\zeta_*$, then
\[
	\zeta<\widetilde P_\alpha(\zeta)<\widetilde P_\alpha(\zeta_*)=\zeta_*.
\]
Therefore $\widetilde P_\alpha^n(\zeta)$ increases to $\zeta_*$. By the $\omega$-periodicity of the coefficients and the uniqueness of solutions, it then follows that
\[
	W(n\omega+\tau;\alpha,\zeta)
	=
	W\bigl(\tau;\alpha,\widetilde P_\alpha^n(\zeta)\bigr),
	\qquad
	\tau\in[0,\omega].
\]
Since $\widetilde P_\alpha^n(\zeta)\to\zeta_*$, the continuous dependence of solutions on initial data, uniformly for $\tau\in[0,\omega]$, yields
\[
	\sup_{\tau\in[0,\omega]}
	\left|
	W(n\omega+\tau;\alpha,\zeta)-W^*(\tau)
	\right|
	\longrightarrow0.
\]
For every $t\ge0$, write $t=n\omega+\tau$, where $\tau\in[0,\omega)$. Since $W^*$ is $\omega$-periodic,
\[
	\left|
	W(t;\alpha,\zeta)-W^*(t)
	\right|
	\le
	\sup_{s\in[0,\omega]}
	\left|
	W(n\omega+s;\alpha,\zeta)-W^*(s)
	\right|.
\]
As $t\to\infty$, we have $n\to\infty$, and hence
\[
	\lim_{t\to\infty}
	\left|
	W(t;\alpha,\zeta)-W^*(t)
	\right|
	=0.
\]
If $\zeta>\zeta_*$, the same argument as in case (i) shows that the iterates eventually leave $\mathcal I_\alpha$; hence the solution blows up in finite forward time. This proves (ii).
	
Finally suppose $\alpha>\alpha_*$. Then $\Theta(\alpha)<0$, and $G_\alpha$ has exactly two zeros $0<\zeta_{\alpha,-}<\zeta_{\alpha,+}$. By strict convexity and the endpoint behavior, we have
\[
	G_\alpha(\zeta)>0
	\quad\text{on }
	(0,\zeta_{\alpha,-})\cup(\zeta_{\alpha,+},R_\alpha),
\]
and
\[
	G_\alpha(\zeta)<0
	\quad\text{on }
	(\zeta_{\alpha,-},\zeta_{\alpha,+}).
\]
At the first zero, $G_\alpha$ crosses from positive to negative, so $G_\alpha'(\zeta_{\alpha,-})<0$. Since $\widetilde P_\alpha'>0$, this gives $0<\widetilde P_\alpha'(\zeta_{\alpha,-})<1$. At the second zero, $G_\alpha$ crosses from negative to positive, so $G_\alpha'(\zeta_{\alpha,+})>0$, and hence $\widetilde P_\alpha'(\zeta_{\alpha,+})>1$. If $0<\zeta<\zeta_{\alpha,-}$, then $\zeta<\widetilde P_\alpha(\zeta)<\zeta_{\alpha,-}$. If $\zeta_{\alpha,-}<\zeta<\zeta_{\alpha,+}$, then $\zeta_{\alpha,-}<\widetilde P_\alpha(\zeta)<\zeta$. Therefore $\widetilde P_\alpha^n(\zeta)\to\zeta_{\alpha,-}$ for every $0<\zeta<\zeta_{\alpha,+}$. Applying the same argument as in the case $\alpha=\alpha_*$, the convergence $\widetilde P_\alpha^n(\zeta)\longrightarrow\zeta_{\alpha,-}$ implies
\[
	\lim_{t\to\infty}
	\left|
	W(t;\alpha,\zeta)-W_{\alpha,-}(t)
	\right|
	=0.
\]
If $\zeta>\zeta_{\alpha,+}$, then $\widetilde P_\alpha(\zeta)>\zeta$, and the same finite-exit argument as in case (i) shows finite-time blow-up.
\end{proof}

\subsection{Turing instability conditions}
\label{pthirdsubsec-turing-instability-conditions}

In this subsection, we restrict our attention to the case $m=2$. By the reduction argument in Section~\ref{pthirdsec-preliminaries}, it suffices to consider the following system:

\begin{equation}\label{pthird-linear-alpha-ODE}
	\frac{\mathrm{d} u}{\mathrm{d} t} = A(t)u-\alpha D(t)u,
\end{equation}
where
\[
u = (x, y), \quad 
A(t) = \begin{bmatrix} 
	a(t) & b(t) \\ 
	c(t) & e(t) 
\end{bmatrix}, \quad 
D(t) = \begin{bmatrix} 
	d_1(t) & 0 \\ 
	0 & d_2(t) 
\end{bmatrix}.
\]
Throughout this and the next subsection,  $a,b,c,e,d_1,d_2\in C_\omega^1(\mathbb R)$, and  $d_i(t)>0$ for all \(t\in\mathbb R\) and \(i=1,2\). Assume, in addition, that \(b\) and \(c\) are nowhere zero and have fixed signs. Then system \eqref{pthird-linear-alpha-ODE} falls into one of the following three types:
\begin{enumerate}
	\item[(i)] cooperative type, that is, $b(t)>0$ and $c(t)>0$ for all $t\in [0,\omega]$;
	\item[(ii)] competitive type, that is, $b(t)<0$ and $c(t)<0$ for all $t\in [0,\omega]$;
	\item[(iii)] predator--prey type, that is, $b(t)c(t)<0$ for all $t\in[0,\omega]$.
\end{enumerate}
By Theorem \ref{pthirdthe-nonlinear-nonoccurrence} and Remark \ref{pthirdrem-general-cone}, Turing instabilities do not occur in cases (i) and (ii). Thus, it remains to consider the predator-prey case. Without loss of generality, we assume that $b(t)>0$ and $c(t)<0$ for all $t\in [0,\omega]$.

We first simplify system \eqref{pthird-linear-alpha-ODE}. Set \(X=x\) and \( Y=r(t)y\), where \(r(t)=\sqrt{-b(t)/c(t)}>0\). Since \(b(t)\) and \(c(t)\) are \(\omega\)-periodic and do not change sign, \(r(t)\) is also a positive \(\omega\)-periodic function. Define \(	q(t)=\sqrt{-b(t)c(t)}>0\). It follows that 
\[
\frac{b(t)}{r(t)}=q(t),
\qquad
r(t)c(t)=-q(t).
\]
Now we compute the transformed system. Since \(X=x\), we have
\[
X'
=
x'
=
\bigl(a(t)-\alpha d_1(t)\bigr)X
+
b(t)\frac{Y}{r(t)}
=
\bigl(a(t)-\alpha d_1(t)\bigr)X
+
q(t)Y.
\]
For the second equation, using \(y=Y/r(t)\), we get
\[
Y'
=
\frac{r'(t)}{r(t)}Y
+
r(t)c(t)X
+
\bigl(e(t)-\alpha d_2(t)\bigr)Y.
\]
Since \(r(t)c(t)=-q(t)\), this becomes
\[
Y'
=
-q(t)X
+
\left(
e(t)-\alpha d_2(t)+\frac{r'(t)}{r(t)}
\right)Y.
\]
Therefore system \eqref{pthird-linear-alpha-ODE} is transformed into
\begin{equation}\label{pthirdsys-skew-symmetric-system}
	\begin{cases}
		X'=p_\alpha(t)X+q(t)Y,\\[1mm]
		Y'=-q(t)X+s_\alpha(t)Y,
	\end{cases}
\end{equation}
where $p_\alpha(t)=a(t)-\alpha d_1(t)$ and $s_\alpha(t)=e(t)-\alpha d_2(t)+\gamma(t)$, with
\[
\gamma(t)=\frac{r'(t)}{r(t)}
=
\frac12
\left(
\frac{b'(t)}{b(t)}
-
\frac{c'(t)}{c(t)}
\right).
\]
Since \(r(t)\) is \(\omega\)-periodic, we have
\[
\int_0^\omega \gamma(t)\,\mathrm dt
=
\int_0^\omega \frac{r'(t)}{r(t)}\,\mathrm dt
=
0.
\]
Thus the transformed coefficient matrix has the form
\begin{equation*}
	A_\alpha(t)
	=
	\begin{pmatrix}
		p_\alpha(t) & q(t)\\
		-q(t) & s_\alpha(t)
	\end{pmatrix}.
\end{equation*}
In particular, the off-diagonal part is skew-symmetric.

Note that the above change of variables is a \(\omega\)-periodic invertible linear transformation. Indeed, if
\[
U(t)=
\begin{pmatrix}
	X(t)\\
	Y(t)
\end{pmatrix},
\qquad
u(t)=
\begin{pmatrix}
	x(t)\\
	y(t)
\end{pmatrix},
\]
then
\[
U(t)=Q(t)u(t),
\qquad
Q(t)=
\begin{pmatrix}
	1&0\\
	0&r(t)
\end{pmatrix}.
\]
Since \(r(t)>0\) and \(r(t+\omega)=r(t)\), we have $Q(t+\omega)=Q(t)$. Let  \(\widetilde T_\alpha\) be the  Poincar\'e map of the transformed system \eqref{pthirdsys-skew-symmetric-system}. Then $\widetilde T_\alpha=Q(0)T_\alpha Q(0)^{-1}$. Hence \(T_\alpha\) and \(\widetilde T_\alpha\) are similar matrices. In particular, they have the same Floquet multipliers, and therefore $\rho(\widetilde T_\alpha)=\rho(T_\alpha)$. Thus the transformation preserves the Floquet stability or instability.

For solutions with \(X(t)\neq 0\), we introduce the slope variable $Z(t)=-Y(t)/X(t)$. We compute the equation for \(Z\). Since \(Y=-ZX\), we have
\[
X'
=
p_\alpha(t)X+q(t)Y
=
\bigl(p_\alpha(t)-q(t)Z\bigr)X.
\]
Also,
\[
Y'
=
-q(t)X+s_\alpha(t)Y
=
-\bigl(q(t)+s_\alpha(t)Z\bigr)X.
\]
On the other hand,
\[
Y'
=
-Z' X-Z X'.
\]
Therefore
\[
-Z' X
-
Z\bigl(p_\alpha(t)-q(t)Z\bigr)X
=
-\bigl(q(t)+s_\alpha(t)Z\bigr)X.
\]
Dividing by \(X\), we obtain
\[
Z'
=
q(t)\bigl(1+Z^2\bigr)-\bigl(p_\alpha(t)-s_\alpha(t)\bigr)Z.
\]
This may also be written as
\begin{equation}\label{pthirdsys-Riccati-alpha}
	Z'
	=
	q(t)\bigl(1+Z^2\bigr)
	-
	\left[
	a(t)-e(t)-\gamma(t)
	+\alpha\bigl(d_2(t)-d_1(t)\bigr)
	\right]Z.
\end{equation}
Moreover, along such a solution, we have 
\[
\frac{X'}{X}
=
p_\alpha(t)-q(t)Z
=
a(t)-\alpha d_1(t)-q(t)Z.
\]

Notice that equation \eqref{pthirdsys-Riccati-alpha} is of the form 
\eqref{pthirdsys-general-Riccati-alpha} with
\[
g(t)=q(t),
\qquad
h(t)=a(t)-e(t)-\gamma(t),
\qquad
\kappa(t)=d_2(t)-d_1(t).
\]
Assume first that $d_2(t)-d_1(t)>0$ for all $t\in[0,\omega]$.
Let $\alpha_*$ denote the critical value given in
Subsection~\ref{pthirdsubsec-periodic-riccati} for this choice of
coefficients. For $\alpha\ge\alpha_*$, define the canonical smaller
positive periodic solution by
\[
Z_\alpha(t)
:=
\begin{cases}
	Z^*(t), & \alpha=\alpha_*,\\[1mm]
	Z_{\alpha,-}(t), & \alpha>\alpha_*.
\end{cases}
\]
Here $Z^*$ and $Z_{\alpha,-}$ are given by
Theorem~\ref{pthirdprop-riccati-global-dynamics}. We also retain the
notation $F_\alpha$, $\widetilde P_\alpha$, $G_\alpha$, and $\mathcal A$ for this
choice of coefficients.

We now use this distinguished positive periodic solution of the Riccati equation to
identify the dominant Floquet multiplier of the modal system. Define
\[
I_\alpha[Z]
:=
\int_0^\omega
\left[p_\alpha(t)-q(t)Z(t)\right]\,\mathrm dt
=
\int_0^\omega
\left[a(t)-\alpha d_1(t)-q(t)Z(t)\right]\,\mathrm dt.
\]
We set $I_-(\alpha):=I_\alpha[Z_\alpha]$. This quantity will be used below
to formulate the modal instability criterion for the case $d_2(t)>d_1(t)$.

We next consider the symmetric Riccati reduction, which will be used when
\(d_1(t)>d_2(t)\). Introduce the rotated variables $\widehat X=Y$ and $\widehat Y=-X$. Then system \eqref{pthirdsys-skew-symmetric-system} becomes
\begin{equation}\label{pthird-skew-symmetric-system-Y}
	\begin{cases}
		{\widehat X}'=s_\alpha(t)\widehat X+q(t)\widehat Y,\\[1mm]
		{\widehat Y}'=-q(t)\widehat X+p_\alpha(t)\widehat Y.
	\end{cases}
\end{equation}
For solutions with \(\widehat X(t)\neq0\), define the rotated slope variable $
\widehat Z(t):=-\frac{\widehat Y(t)}{\widehat X(t)}$. Then \(\widehat Z\) satisfies the Riccati equation
\[
\widehat Z'
=
q(t)\bigl(1+\widehat Z^2\bigr)
-
\bigl[s_\alpha(t)-p_\alpha(t)\bigr]\widehat Z .
\]
Equivalently,

\begin{equation}\label{pthirdsys-Riccati-alpha-hat}
	\widehat Z'
	=
	q(t)\bigl(1+\widehat Z^2\bigr)
	-
	\left[
	e(t)-a(t)+\gamma(t)
	+
	\alpha\bigl(d_1(t)-d_2(t)\bigr)
	\right]\widehat Z .
\end{equation}
Moreover, we have
\[
\frac{\widehat X'}{\widehat X}
=
s_\alpha(t)-q(t)\widehat Z .
\]
Thus the transformed system has the same form as \eqref{pthirdsys-skew-symmetric-system}, with \(p_\alpha\) and \(s_\alpha\) interchanged. Define
\[
\widehat I_\alpha[Z]
:=
\int_0^\omega
\left[
s_\alpha(t)-q(t)Z(t)
\right]\,\mathrm dt
=
\int_0^\omega
\left[
e(t)-\alpha d_2(t)-q(t)Z(t)
\right]\,\mathrm dt.
\]

Assume now that $d_1(t)-d_2(t)>0$ for all $t\in[0,\omega]$.
Equation \eqref{pthirdsys-Riccati-alpha-hat} is again of the form
\eqref{pthirdsys-general-Riccati-alpha} with
\[
g(t)=q(t),
\qquad
h(t)=e(t)-a(t)+\gamma(t),
\qquad
\kappa(t)=d_1(t)-d_2(t).
\]
Applying Subsection~\ref{pthirdsubsec-periodic-riccati} with these
coefficients, we denote the corresponding critical value by
$\widehat\alpha_*$ and define $\widehat Z_\alpha$ analogously for
$\alpha\ge\widehat\alpha_*$. Set
\[
\widehat I_-(\alpha)
:=\widehat I_\alpha[\widehat Z_\alpha].
\]

\begin{lemma}\label{pthirdlem-postcritical-exact-modal-criterion}
The following statements are valid:
\begin{enumerate}
	\item[(i)] If   \(d_2(t)-d_1(t)>0\) for all \(t\in[0,\omega]\), then for every \(\alpha\ge\alpha_*\), we have $\rho(T_\alpha)=e^{I_-(\alpha)}$. In particular, $\rho(T_\alpha)>1$ if and only if $I_-(\alpha)>0$.
		
	\item[(ii)] If \(d_1(t)-d_2(t)>0\) for all \(t\in[0,\omega]\), then for every \(\alpha\ge\widehat\alpha_*\), we have $\rho(T_\alpha)=e^{\widehat I_-(\alpha)}$. In particular, $\rho(T_\alpha)>1$ if and only if $\widehat I_-(\alpha)>0$.
\end{enumerate}
\end{lemma}

\begin{proof}
We first prove (i). Since \(T_\alpha\) and \(\widetilde T_\alpha\) are similar, it is enough to work with the transformed system \eqref{pthirdsys-skew-symmetric-system}.
	
Suppose that \(\alpha>\alpha_*\). By Theorem~\ref{pthirdprop-riccati-global-dynamics}, equation \eqref{pthirdsys-Riccati-alpha} has exactly two positive \(\omega\)-periodic solutions \(Z_{\alpha,-}\) and \(Z_{\alpha,+}\) with $Z_{\alpha,-}(t)<Z_{\alpha,+}(t)$ for all $t\in\mathbb R$. For \(\iota\in\{-,+\}\), define
\[
	X_{\alpha,\iota}(t)
	:=
	\exp\left\{
	\int_0^t
	\bigl[p_\alpha(\tau)-q(\tau)Z_{\alpha,\iota}(\tau)\bigr]\,\mathrm d\tau
	\right\},
	\qquad
	Y_{\alpha,\iota}(t):=-Z_{\alpha,\iota}(t)X_{\alpha,\iota}(t).
\]
Then
\[
	\frac{X_{\alpha,\iota}'}{X_{\alpha,\iota}}
	=
	p_\alpha(t)-q(t)Z_{\alpha,\iota}(t).
\]
Using the Riccati equation satisfied by \(Z_{\alpha,\iota}\), we obtain
\[
	Y_{\alpha,\iota}'
	=
	-Z_{\alpha,\iota}'X_{\alpha,\iota}
	-
	Z_{\alpha,\iota}X_{\alpha,\iota}'
	=
	-q(t)X_{\alpha,\iota}+s_\alpha(t)Y_{\alpha,\iota}.
\]
Thus \((X_{\alpha,\iota},Y_{\alpha,\iota})\) solves equation \eqref{pthirdsys-skew-symmetric-system}. Since \(Z_{\alpha,\iota}\) is \(\omega\)-periodic,
\[
	\begin{pmatrix}
		X_{\alpha,\iota}(\omega)\\
		Y_{\alpha,\iota}(\omega)
	\end{pmatrix}
	=
	e^{I_\alpha[Z_{\alpha,\iota}]}
	\begin{pmatrix}
		X_{\alpha,\iota}(0)\\
		Y_{\alpha,\iota}(0)
	\end{pmatrix}.
\]
Hence $e^{I_\alpha[Z_{\alpha,-}]}$ and $e^{I_\alpha[Z_{\alpha,+}]}$ are Floquet multipliers of the transformed system. The corresponding initial vectors $(1,-Z_{\alpha,-}(0))$ and $(1,-Z_{\alpha,+}(0))$ are linearly independent since \(Z_{\alpha,-}(0)\neq Z_{\alpha,+}(0)\). Hence these are exactly the two Floquet multipliers of \(\widetilde T_\alpha\). Moreover,
\[
	I_\alpha[Z_{\alpha,-}]
	-
	I_\alpha[Z_{\alpha,+}]
	=
	\int_0^\omega
	q(t)\bigl[
	Z_{\alpha,+}(t)-Z_{\alpha,-}(t)
	\bigr]\,\mathrm dt
	>0.
\]
Therefore
\[
	\rho(\widetilde T_\alpha)
	=
	e^{I_\alpha[Z_{\alpha,-}]}
	=
	e^{I_-(\alpha)}.
\]
Since \(T_\alpha\) and \(\widetilde T_\alpha\) are similar, we obtain $\rho(T_\alpha)=e^{I_-(\alpha)}$.
	
It remains to consider the critical case \(\alpha=\alpha_*\). Let \(Z^*\) be the unique positive \(\omega\)-periodic solution, and set $\mu:=e^{I_{\alpha_*}[Z^*]}=e^{I_-(\alpha_*)}$. The same construction as above shows that \(\mu\) is one Floquet multiplier of \(\widetilde T_{\alpha_*}\).	We claim that the other Floquet multiplier is also \(\mu\). Indeed, the Riccati variational equation gives
\[
	P_{\alpha_*}'(Z^*(0))
	=
	\exp\left\{
	\int_0^\omega
	\partial_ZF_{\alpha_*}(t,Z^*(t))\,\mathrm dt
	\right\}.
\]
At the critical value \(\alpha=\alpha_*\), the function \(G_{\alpha_*}\) has a unique zero at an interior minimum. Hence $G_{\alpha_*}'(Z^*(0))=0$ and $P_{\alpha_*}'(Z^*(0))=1$. On the other hand,
\[
	\partial_ZF_{\alpha_*}(t,Z^*(t))
	=
	2q(t)Z^*(t)-\bigl[p_{\alpha_*}(t)-s_{\alpha_*}(t)\bigr].
\]
By Liouville's formula, it follows that 
\[
	\det \widetilde T_{\alpha_*}
	=
	\exp\left\{
	\int_0^\omega
	\bigl[p_{\alpha_*}(t)+s_{\alpha_*}(t)\bigr]\,\mathrm dt
	\right\},
\]
and since
\[
     \mu
	=
	\exp\left\{
	\int_0^\omega
	\bigl[p_{\alpha_*}(t)-q(t)Z^*(t)\bigr]\,\mathrm dt
	\right\},
\]
we obtain
\[
	P_{\alpha_*}'(Z^*(0))
	=
	\frac{\det \widetilde T_{\alpha_*}}{\mu^2}.
\]
Since \(P_{\alpha_*}'(Z^*(0))=1\), it follows that $\det \widetilde T_{\alpha_*}=\mu^2$. The product of the two Floquet multipliers is \(\det \widetilde T_{\alpha_*}\). Since one multiplier is \(\mu\), the other one is also \(\mu\). Thus $\rho(\widetilde T_{\alpha_*})=\mu=e^{I_-(\alpha_*)}$. Again, by similarity, $\rho(T_{\alpha_*})=e^{I_-(\alpha_*)}$. Combining the cases \(\alpha>\alpha_*\) and \(\alpha=\alpha_*\), we get $\rho(T_\alpha)=e^{I_-(\alpha)}$ for every $\alpha\ge\alpha_*$. Therefore $\rho(T_\alpha)>1$ if and only if $I_-(\alpha)>0$. The proof of (ii) is the same after applying the preceding argument to the rotated system \eqref{pthird-skew-symmetric-system-Y}. 
\end{proof}

We next establish the dependence of the canonical solution of the Riccati equation
and the associated growth integral on the modal parameter.

\begin{lemma}\label{pthirdlem-postcritical-concavity}
Let  $d_2(t)-d_1(t)>0$ for all $t\in[0,\omega]$. Then $I_-$ is continuous on $[\alpha_*,\infty)$, strictly concave on $(\alpha_*,\infty)$, and $\lim_{\alpha\to\infty}I_-(\alpha)=-\infty$. The analogous conclusions hold for $\widehat I_-$ when $d_1(t)-d_2(t)>0$.
\end{lemma}

\begin{proof}
We prove the assertions for \(d_2(t)>d_1(t)\). In this case, $\kappa(t)=d_2(t)-d_1(t)$. For each \(\alpha>\alpha_*\), Theorem~\ref{pthirdprop-riccati-global-dynamics} gives $0<P_\alpha'(Z_\alpha(0))<1$. Hence,
\[
	\partial_\zeta\bigl(P_\alpha(\zeta)-\zeta\bigr)
	\big|_{\zeta=Z_\alpha(0)}
	=
	P_\alpha'(Z_\alpha(0))-1\neq0.
\]
By the implicit function theorem, together with the smooth dependence of solutions of the Riccati equation on the parameter and the initial value, \(Z_\alpha\) depends \(C^2\)-smoothly on \(\alpha\in(\alpha_*,\infty)\). Since 
\[
	\partial_ZF_\alpha(t,Z_\alpha(t))
	=
	2q(t)Z_\alpha(t)
	-
	\left[
	a(t)-e(t)-\gamma(t)+\alpha\kappa(t)
	\right],
\]
the variational equation for the period map yields
\[
	P_\alpha'(Z_\alpha(0))
	=
	\exp\left\{
	\int_0^\omega
	\partial_ZF_\alpha(t,Z_\alpha(t))\,\mathrm dt
	\right\}.
\]
Thus,
\[
	\int_0^\omega
	\partial_ZF_\alpha(t,Z_\alpha(t))\,\mathrm dt
	=
	\log P_\alpha'(Z_\alpha(0))
	<0.
\]
Let \(w_\alpha:=\partial_\alpha Z_\alpha\). Differentiating the Riccati equation with respect to \(\alpha\), we obtain
\[
	w_\alpha'
	=
	\partial_ZF_\alpha(t,Z_\alpha(t))w_\alpha
	-
	\kappa(t)Z_\alpha(t).
\]
Since \(\partial_ZF_\alpha(t,Z_\alpha(t))\) is \(\omega\)-periodic and has negative mean, this equation admits a unique \(\omega\)-periodic solution given by
\[
	w_\alpha(t)
	=
	-
	\int_{-\infty}^t
	\exp\left\{
	\int_s^t
	\partial_ZF_\alpha(\tau,Z_\alpha(\tau))\,\mathrm d\tau
	\right\}
	\kappa(s)Z_\alpha(s)\,\mathrm ds
	<0.
\]
Let $v_\alpha:=\partial_\alpha^2 Z_\alpha$. Differentiating once more gives
\[
	v_\alpha'
	=
	\partial_ZF_\alpha(t,Z_\alpha(t))v_\alpha
	+
	2q(t)w^2_\alpha(t)
	-
	2\kappa(t)w_\alpha(t).
\]
Since \(w_\alpha<0\) and \(\kappa>0\), the forcing term satisfies
\[
	2q(t)w^2_\alpha(t)
	-
	2\kappa(t)w_\alpha(t)>0.
\]
Hence the unique periodic solution satisfies
\[
	v_\alpha(t)
	=
	\int_{-\infty}^t
	\exp\left\{
	\int_s^t
	\partial_ZF_\alpha(\tau,Z_\alpha(\tau))\,\mathrm d\tau
	\right\}
	\left[
	2q(s)w^2_\alpha(s)
	-
	2\kappa(s)w_\alpha(s)
	\right]
	\,\mathrm ds
	>0.
\]
It follows that
\[
	I_-''(\alpha)
	=
	-
	\int_0^\omega q(t)v_\alpha(t)\,\mathrm dt<0.
\]
So \(I_-\) is strictly concave on \((\alpha_*,\infty)\).

By Lemma~\ref{pthirdlem-postcritical-exact-modal-criterion}, we have $I_-(\alpha)=\log\rho(T_\alpha)$ for all $\alpha\ge\alpha_*$. Since the monodromy matrix depends continuously on \(\alpha\), so does its spectral radius. Thus \(I_-\) is continuous at \(\alpha_*\), and hence on \([\alpha_*,\infty)\). Finally, since \(q(t)Z_\alpha(t)>0\), 
\[ 
	I_-(\alpha) < \int_0^\omega a(t)\,\mathrm dt - \alpha\int_0^\omega d_1(t)\,\mathrm dt \longrightarrow-\infty. 
\]
The proof for the rotated Riccati equation is analogous, with the roles of \(d_1\) and \(d_2\) interchanged.
\end{proof}

\begin{theorem}\label{pthirdthm-refined-turing-instability}
Assume that (A1) and (A2) hold, and that one of the following two alternatives holds:
\begin{enumerate}
		\item[(i)] \(d_2(t)-d_1(t)>0\) for all \(t\in[0,\omega]\), and there exists \(\alpha>0\) such that $\alpha\ge\alpha_*$ and $I_-(\alpha)>0$.
		
		\item[(ii)] \(d_1(t)-d_2(t)>0\) for all \(t\in[0,\omega]\), and there exists \(\alpha>0\) such that $\alpha\ge\widehat\alpha_*$ and $\widehat I_-(\alpha)>0$.
\end{enumerate}
Then a Turing instability occurs on a suitable rescaling of the spatial domain \(\Omega\). Furthermore, if the above \(\alpha\) can be chosen as a nonzero Neumann eigenvalue \(\lambda_j\), for some \(j\ge1\), then a Turing instability occurs on the original domain \(\Omega\).
\end{theorem}

\begin{proof}
By assumption (A2), the spatially homogeneous periodic solution \(u^*(t)\) is stable as a solution of the kinetic system~\eqref{pthirdsys-nonlinear-ODE}. In either case, Lemma~\ref{pthirdlem-postcritical-exact-modal-criterion} gives \(\rho(T_\alpha)>1\), since \(\rho(T_\alpha)=e^{I_-(\alpha)}\) in case (i) and
\(\rho(T_\alpha)=e^{\widehat I_-(\alpha)}\) in case (ii). Thus the modal equation corresponding to this value of \(\alpha\) is unstable.

It remains only to realize this value of \(\alpha\) as a nonzero Neumann eigenvalue after rescaling the spatial domain. Let \(\lambda_j>0\) be any nonzero Neumann eigenvalue of the original domain \(\Omega\). Under the
rescaling $\widetilde\Omega:=\beta\Omega$, the corresponding Neumann eigenvalue becomes $\widetilde\lambda_j=\beta^{-2}\lambda_j$. Choosing $\beta=\sqrt{\frac{\lambda_j}{\alpha}}$ gives \(\widetilde\lambda_j=\alpha\). Therefore the unstable modal equation corresponds to a genuine spatially nonhomogeneous perturbation of the reaction-diffusion system on \(\widetilde\Omega\). Since the homogeneous mode is stable by (A2), while this nonzero spatial mode is unstable, a Turing instability occurs on the rescaled domain. If, in addition, this \(\alpha\) is a nonzero Neumann eigenvalue \(\lambda_j\) on \(\Omega\), then the unstable nonhomogeneous mode is already available on the original domain. Hence a Turing instability occurs on \(\Omega\).
\end{proof}

\begin{remark}\label{pthirdrem-supersolution}
In applications, it is not necessary to determine the positive periodic solution of the Riccati equation explicitly. Under the diffusion ordering in Theorem~\ref{pthirdthm-refined-turing-instability}(i), let $\alpha>0$ and let $\overline Z$ be a positive $\omega$-periodic supersolution of \eqref{pthirdsys-Riccati-alpha} such that $I_\alpha[\overline Z]>0$. The same instability conclusion remains valid. Indeed, a straightforward comparison gives
\[
	P_\alpha\bigl([0,\overline Z(0)]\bigr)
	\subset[0,\overline Z(0)],
	\qquad P_\alpha(0)>0.
\]
Thus,  there exists a positive periodic solution below $\overline Z$. By Theorem~\ref{pthirdprop-riccati-global-dynamics}, we have $\alpha\ge\alpha_*$ and $Z_\alpha(t)\le\overline Z(t)$ for all $t\in\mathbb R$. Consequently,
$I_-(\alpha)=I_\alpha[Z_\alpha]\ge I_\alpha[\overline Z]>0$. The same argument applies to the rotated Riccati equation in case (ii).
\end{remark}

\begin{remark}
The upper bound in Lemma~\ref{pthirdlem-alpha-set-structure},  applied with the corresponding coefficients $g,h,\kappa$, gives explicit upper bounds for $\alpha_*$ and $\widehat\alpha_*$. These provide a convenient way to check the criterion numerically. In practice, we may choose $\alpha>0$ above the corresponding explicit upper bound and then compute $I_-(\alpha)$ or $\widehat I_-(\alpha)$. A positive value guarantees Turing instability by Theorem~\ref{pthirdthm-refined-turing-instability}.
\end{remark}

\begin{remark}\label{pthirdrem-Van-Gorder}
Van Gorder \cite[Theorem~2.1]{VanGorder-PRSA-2020} proposed a sufficient condition for instability in two-dimensional nonautonomous reaction--diffusion systems. In the periodic setting, the corresponding conditions reduce to that there exists  a modal parameter $\alpha>0$ such that either
\[
	\det(A(t)-\alpha D(t))
	<
	b(t)\frac{\mathrm d}{\mathrm dt}
	\left(
	\frac{a(t)-\alpha d_1(t)}{b(t)}
	\right)
	\quad\text{for all }t\in[0,\omega],
\]
or
\[
	\det(A(t)-\alpha D(t))
	<
	c(t)\frac{\mathrm d}{\mathrm dt}
	\left(
	\frac{e(t)-\alpha d_2(t)}{c(t)}
	\right)
	\quad\text{for all }t\in[0,\omega].
\]
Within our predator--prey setting  here, these pointwise conditions give particular strict periodic supersolution conditions for one of the two associated Riccati equations, and hence, are covered by our Riccati type criterion. Our criterion, however, only requires a positive growth integral over one period. This distinction is illustrated in section~\ref{pthirdexp-segel-jackson}, where our criterion detects an unstable mode, while both pointwise conditions fail at \(t=\omega/2\) for sufficiently small \(\varepsilon>0\).
\end{remark}

\subsection{Revisiting autonomous systems}
We now compare the criterion obtained for periodic systems with two autonomous situations, according to whether $u^*$ is a constant positive equilibrium or a nonconstant periodic orbit of system \eqref{pthirdsys-nonlinear-ODE}. In the first case, the modal stability problem reduces to an algebraic eigenvalue problem. We will show below that our criterion then becomes exactly the classical  Turing condition. In the second case, we first establish the modal structure for periodic orbits of autonomous 
systems and then compare our criterion with a criterion for  Hopf-Turing instability obtained in  \cite{Yi-JDE-2021}. We also explain that these approaches are complementary.

We first consider the case of an equilibrium for an autonomous system. Then
\[
A(t)\equiv A=
\begin{pmatrix}
	a&b\\
	c&e
\end{pmatrix},
\qquad
D(t)\equiv D=
\begin{pmatrix}
	d_1&0\\
	0&d_2
\end{pmatrix},
\qquad d_1,d_2>0 .
\]
We retain the predator--prey assumption $bc<0$ and assume that the kinetic equilibrium is linearly asymptotically stable, namely
\begin{equation}\label{pthird-autonomous-kinetic-stability}
	\operatorname{tr}A=a+e<0,
	\qquad
	\det A=ae-bc>0 .
\end{equation}
For a spatial mode with parameter \(\alpha>0\), the modal matrix is
\[
A-\alpha D
=
\begin{pmatrix}
	a-\alpha d_1&b\\
	c&e-\alpha d_2
\end{pmatrix}.
\]
Since
\[
\operatorname{tr}(A-\alpha D)
=
a+e-\alpha(d_1+d_2)<0,
\]
the mode is linearly unstable if and only if
\begin{equation}\label{pthird-autonomous-modal-condition}
	H(\alpha):=\det(A-\alpha D)<0 .
\end{equation}
Explicitly,
\begin{equation}\label{pthird-autonomous-H-alpha}
	H(\alpha)
	=
	d_1d_2\alpha^2
	-
	(d_2a+d_1e)\alpha
	+
	\det A .
\end{equation}
Therefore, in the present notation, the classical Turing
criterion is
\begin{equation}\label{pthird-autonomous-continuous-classical}
	d_2a+d_1e>0,
	\qquad
	(d_2a+d_1e)^2-4d_1d_2\det A>0.
\end{equation}
The unstable interval is $\alpha_-<\alpha<\alpha_+$,
where
\[
\alpha_\pm
=
\frac{
	d_2a+d_1e
	\pm
	\sqrt{(d_2a+d_1e)^2-4d_1d_2\det A}
}{
	2d_1d_2
}.
\]
This is the usual continuous version formulation of the classical two-species Turing criterion \cite[Section~2.3]{Murray-book-2003}. On a fixed spatial domain, the formulation is
\begin{equation}\label{pthird-autonomous-fixed-domain-classical}
	H(\lambda_j)=\det(A-\lambda_jD)<0
	\qquad
	\text{for some }j\ge1.
\end{equation}
Equivalently, $\lambda_j\in(\alpha_-,\alpha_+)$ for some $j\ge1$ \cite[Theorem~1]{Casten-Holland-SIAP-1977}. The first asks whether the unstable interval is nonempty, while the second asks whether the spectrum of the given domain actually intersects this interval.

We now show that our criterion reduces exactly to the determinant condition
\eqref{pthird-autonomous-modal-condition}. Set
\[
q:=\sqrt{-bc}>0,
\qquad
p_\alpha:=a-\alpha d_1,
\qquad
s_\alpha:=e-\alpha d_2 .
\]
Since \(b\) and \(c\) are constant, the transformation used above is
constant. The transformed modal system is
\begin{equation}\label{pthird-autonomous-transformed-system}
	\begin{cases}
		X'=p_\alpha X+qY,\\[1mm]
		Y'=-qX+s_\alpha Y.
	\end{cases}
\end{equation}
Its coefficient matrix
\[
A_\alpha
=
\begin{pmatrix}
	p_\alpha&q\\
	-q&s_\alpha
\end{pmatrix}
\]
is similar to \(A-\alpha D\). Hence \(A_\alpha\) and \(A-\alpha D\) have
the same eigenvalues.

If \(d_2>d_1\), the Riccati equation associated with
\eqref{pthird-autonomous-transformed-system} is
\begin{equation}\label{pthird-autonomous-Riccati}
	Z'
	=
	q(1+Z^2)
	-
	(p_\alpha-s_\alpha)Z .
\end{equation}
Every periodic solution of this autonomous scalar equation is constant.
The critical value defined earlier therefore reduces to
\[
	\alpha_*
	=
	\frac{2q-a+e}{d_2-d_1}.
\]
For \(\alpha\ge\alpha_*\), the smaller positive equilibrium is
\[
Z_\alpha
=
\frac{
	p_\alpha-s_\alpha
	-
	\sqrt{(p_\alpha-s_\alpha)^2-4q^2}
}{
	2q
}.
\]
Similarly, if \(d_1>d_2\), then the critical value for the rotated
Riccati equation reduces to
\begin{equation}\label{pthird-autonomous-hat-alpha-star}
	\widehat\alpha_*
	=
	\frac{2q+a-e}{d_1-d_2}.
\end{equation}

\begin{proposition}\label{pthirdprop-autonomous-reduction}
	Assume that
	\[
	\operatorname{tr}A<0,
	\qquad
	\det A>0,
	\qquad
	bc<0 .
	\]
	Then the following statements are valid:
	
	\begin{enumerate}
		\item[(i)]
		If \(d_2>d_1\), then for every \(\alpha>0\),
		\[
		\det(A-\alpha D)<0
		\quad\Longleftrightarrow\quad
		\alpha\ge\alpha_*
		\ \text{and}\
		I_-(\alpha)>0 .
		\]
		
		\item[(ii)]
		If \(d_1>d_2\), then for every \(\alpha>0\),
		\[
		\det(A-\alpha D)<0
		\quad\Longleftrightarrow\quad
		\alpha\ge\widehat\alpha_*
		\ \text{and}\
		\widehat I_-(\alpha)>0 .
		\]
	\end{enumerate}
\end{proposition}

\begin{proof}
Fix \(\alpha>0\) and first suppose that \(d_2>d_1\). If \(\alpha\ge\alpha_*\), then \(Z_\alpha\) is the smaller
positive root of
\[
	q(1+Z^2)-(p_\alpha-s_\alpha)Z=0.
\]
Therefore
\[
	A_\alpha
	\begin{pmatrix}
		1\\
		-Z_\alpha
	\end{pmatrix}
	=
	\bigl(p_\alpha-qZ_\alpha\bigr)
	\begin{pmatrix}
		1\\
		-Z_\alpha
	\end{pmatrix}.
\]
Since \(Z_\alpha\) is the smaller positive root, the corresponding eigenvalue is the larger eigenvalue of \(A_\alpha\). Denote this eigenvalue by \(\lambda_+(\alpha)\). Then $\lambda_+(\alpha)=p_\alpha-qZ_\alpha$. By the definition of \(I_-(\alpha)\), it follows that
\[
	I_-(\alpha)
	=
	\int_0^\omega
	\bigl(p_\alpha-qZ_\alpha\bigr)\,\mathrm dt
	=
	\omega\lambda_+(\alpha).
\]
Since \(A_\alpha\) is similar to \(A-\alpha D\), \(\lambda_+(\alpha)\) is also the larger eigenvalue of \(A-\alpha D\). Moreover, $\operatorname{tr}(A-\alpha D)<0$. Hence
\[
	I_-(\alpha)>0
	\quad\Longleftrightarrow\quad
	\lambda_+(\alpha)>0
	\quad\Longleftrightarrow\quad
	\det(A-\alpha D)<0 .
\]
Thus, under the assumption \(\alpha\ge\alpha_*\), the two conditions
are equivalent.
	
It remains to show that \(\det(A-\alpha D)<0\)  implies
\(\alpha\ge\alpha_*\). Suppose that $H(\alpha)=\det(A-\alpha D)<0$. Since \(H(0)=\det A>0\) and the leading coefficient of \(H\) is positive, this is possible only if $d_2a+d_1e>0$. Using \(a+e<0\), we write
\[
	d_2a+d_1e
	=
	d_1(a+e)+(d_2-d_1)a.
\]
Since \(d_1(a+e)<0\) and \(d_2>d_1\), the above positivity implies \(a>0\). Hence \(e<0\). Therefore $s_\alpha=e-\alpha d_2<0$. Now
\[
	\det A_\alpha
	=
	p_\alpha s_\alpha+q^2
	=
	\det(A-\alpha D)<0.
\]
Thus we get \(p_\alpha>0\). Consequently,
\[
	p_\alpha-s_\alpha
	\ge
	2\sqrt{-p_\alpha s_\alpha}
	>
	2q.
\]
This gives \(\alpha>\alpha_*\), and the equivalence established above
completes the proof of \({\rm(i)}\).

The proof of \({\rm(ii)}\) is the same by using the rotated Riccati equation.
\end{proof}

Proposition~\ref{pthirdprop-autonomous-reduction} shows that, in the
autonomous equilibrium case, the Riccati conditions in
Theorem~\ref{pthirdthm-refined-turing-instability} are precisely
equivalent to the classical modal instability condition $\det(A-\alpha D)<0$. Thus, the periodic criterion recovers exactly the classical determinant criterion. The two conclusions of Theorem~\ref{pthirdthm-refined-turing-instability} correspond, respectively, to the continuous-mode and fixed-domain formulations.

Next, we  consider the case of an periodic orbit  for an autonomous system.
Assume that \(u^*\) is a nonconstant orbitally stable periodic solution of $u'=F(u)$, with least period
\(\omega\). For the reaction--diffusion system with constant diffusion matrix \(D\), the modal equations parametrized by \(\alpha\ge0\) are
\begin{equation}
	\label{pthirdsys-linear-autonomous-periodic-orbit-case}
	v'
	=
	\bigl[A(t)-\alpha D\bigr]v,
	\qquad
	A(t):=D_uF(u^*(t)).
\end{equation}
Thus, the modal problem is still a periodic Floquet problem. As before, we may assume that \(b(t)>0\) and \(c(t)<0\). At \(\alpha=0\), one Floquet multiplier is \(1\). Denote the other multiplier by \(\mu_2(0)\). Orbital stability and Liouville's formula imply $0<\mu_2(0)\le1$.

We first present a consequence of the results on the Riccati equation, which 
applies to a given orbitally stable periodic solution of the
autonomous kinetic system \(u'=F(u)\).

\begin{proposition}\label{pthirdprop-autonomous-periodic-modal-structure}
Suppose that \(d_2>d_1\). Then the following statements are valid:
\begin{enumerate}
	\item[(i)] The Riccati critical value satisfies \(\alpha_*>0\).
		
	\item[(ii)] The modal system is stable at \(\alpha=\alpha_*\). More precisely,
	\begin{equation}\label{pthird-critical-mode-stability}
			\rho(T_{\alpha_*})
			=
			\sqrt{\mu_2(0)}
			\exp\left\{
			-\frac12\alpha_*(d_1+d_2)\omega
			\right\}
			<1.
	\end{equation}
	\item[(iii)] The instability set $\left\{\alpha\ge\alpha_*:\rho(T_\alpha)>1\right\}$ is either empty or a single bounded open interval.
\end{enumerate}
The analogous conclusions hold when \(d_1>d_2\), with \(\widehat\alpha_*\) in place of \(\alpha_*\).
\end{proposition}

\begin{proof}
We first prove that the Riccati equation has no positive periodic solution at \(\alpha=0\). Autonomy implies that $\dot u^*(t)$ is a nontrivial periodic solution of \(v'=A(t)v\). Moreover, we have  \(\dot u^*_1\not\equiv0\). Otherwise, \(b(t)>0\) would imply \(\dot u^*_2\equiv0\). Since \(\dot u^*_1\) is the derivative of a periodic function, it has a zero \(t_0\), and periodicity gives another zero at \(t_0+\omega\). Suppose, to the contrary, that \(Z\) is a positive periodic solution of \eqref{pthirdsys-Riccati-alpha} at \(\alpha=0\). The Riccati reconstruction gives a solution \(v(t)=(v_1(t),v_2(t))\) of \(v'=A(t)v\), where
\[
	v_1(t)
	=
	\exp\left\{
	\int_0^t
	[a(s)-q(s)Z(s)]\,\mathrm ds
	\right\}>0,
	\qquad
	v_2(t)=-\frac{Z(t)v_1(t)}{r(t)}.
\]
Set
\[
	W(t):=v_1(t)\dot u^*_2(t)-v_2(t)\dot u^*_1(t).
\]
Since \(v\) and \(\dot u^*\) solve the same linear system, it follows that
\[
	W'=(a(t)+e(t))W,
	\qquad
	\left(\frac{\dot u^*_1}{v_1}\right)'
	=
	\frac{b(t)W(t)}{v_1(t)^2}.
\]
If \(W\not\equiv0\), then \(W(t)\) is either strictly positive or strictly negative for all \(t\). Thus \(\dot u^*_1/v_1\) is strictly monotone, which is impossible because it vanishes at both \(t_0\) and \(t_0+\omega\). If \(W\equiv0\), then \(v\) and \(\dot u^*\) are linearly dependent, which is also impossible since \(v_1\) is everywhere positive whereas \(\dot u^*_1\) has a zero. Therefore \(0\notin\mathcal A\). Since $\mathcal A=[\alpha_*,\infty)$, by Lemma~\ref{pthirdlem-alpha-set-structure}, we obtain \(\alpha_*>0\).
	
For every \(\alpha\ge0\), Liouville's formula gives
\[
	\det T_\alpha
	=
	\det T_0\,e^{-\alpha(d_1+d_2)\omega}
	=
	\mu_2(0)e^{-\alpha(d_1+d_2)\omega}.
\]
As established in the proof of Lemma~\ref{pthirdlem-postcritical-exact-modal-criterion}, the two Floquet multipliers coincide and are positive at \(\alpha=\alpha_*\). This yields \eqref{pthird-critical-mode-stability}. Finally, Lemmas \ref{pthirdlem-postcritical-exact-modal-criterion} and \ref{pthirdlem-postcritical-concavity} give 
\[
\left\{ \alpha\ge\alpha_*:\rho(T_\alpha)>1\right\} =\left\{\alpha\ge\alpha_*:I_-(\alpha)>0\right\}, 
\]
where \(I_-\) is continuous on \([\alpha_*,\infty)\), strictly concave on \((\alpha_*,\infty)\), and satisfies $I_-(\alpha_*)<0$ and $\lim_{\alpha\to\infty}I_-(\alpha)=-\infty$. Therefore, the set on which \(I_-\) is positive is either empty or a single bounded open interval whose closure lies in \((\alpha_*,\infty)\).

For \(d_1>d_2\), the same argument using the second component of \(\dot u^*\) excludes a positive periodic solution of the rotated Riccati equation at \(\alpha=0\). The remaining conclusions follow from the corresponding rotated results.
\end{proof}

Proposition~\ref{pthirdprop-autonomous-periodic-modal-structure}
describes the entire instability set beyond the critical value
\(\alpha_*\), but does not guarantee that this set is nonempty.
Its nonemptiness can be verified by checking whether
\(I_-(\alpha)>0\) for some \(\alpha\ge\alpha_*\), or by using the
supersolution criterion in Remark~\ref{pthirdrem-supersolution}.

To compare our criterion with the criterion obtained in \cite{Yi-JDE-2021}, we focus on stable periodic orbits generated by a Hopf bifurcation. Consider a one-parameter family of planar kinetic systems
\begin{equation}\label{pthirdsys-Hopf-ODE}
	u'=F(\varepsilon,u),
\end{equation}
where \(F\) is sufficiently smooth. Assume that a nondegenerate Hopf bifurcation occurs at \((\bar\varepsilon,\bar u)\) and generates a branch of
orbitally asymptotically stable periodic solutions. Write
\[
\bar J:=D_uF(\bar\varepsilon,\bar u)
=
\begin{pmatrix}
	\bar a& \bar b\\
	\bar c&-\bar a
\end{pmatrix},
\qquad
\Sigma:=\det \bar J>0.
\]
The eigenvalues of \(\bar J\) are $\pm i\sqrt{\Sigma}$. Consequently, every nontrivial solution of the linearized system $v'=\bar Jv$ is periodic, with least period $\bar\omega:=\frac{2\pi}{\sqrt{\Sigma}}$. Moreover, if \(u_\varepsilon^*\) denotes the periodic solutions generated by
the Hopf bifurcation and \(\omega_\varepsilon\) denotes their least periods,
then $\omega_\varepsilon\to\bar\omega$ as $\varepsilon\to\bar\varepsilon$ along the bifurcating branch. For \(\alpha\geq0\), let \(T_{\varepsilon,\alpha}\) denote the
monodromy matrix over one period \(\omega_\varepsilon\) of
\[
v'
=
\bigl[D_uF(\varepsilon,u_\varepsilon^*(t))-\alpha D\bigr]v.
\]
Note that $\bar b \bar c=-\bar a^2-\Sigma<0$. Hence the Jacobian at a Hopf point is automatically of predator--prey type. In particular, the same sign condition holds along all sufficiently small bifurcating periodic orbits. Set
\[
\bar q:=\sqrt{-\bar b\bar c}
=\sqrt{\bar a^2+\Sigma}>|\bar a|.
\]

Let \(l_1\) be the complex cubic coefficient in the Hopf normal
form. In the absence of cross-diffusion,
condition~(3.24) in  \cite{Yi-JDE-2021} reduces to
\begin{equation}\label{pthird-Yi-diagonal-Hopf-condition}
	\sqrt{\Sigma}\,(d_1+d_2)
	+
	\frac{\operatorname{Im}l_1}{\operatorname{Re}l_1}
	(d_2-d_1)\bar a
	<0.
\end{equation}
For each periodic orbit on this stable branch sufficiently close to
the Hopf point, this condition implies that the Floquet multiplier issuing from the multiplier \(1\) at \(\alpha=0\)
moves outside the unit circle for all sufficiently small positive
\(\alpha\). Consequently, the spatially homogeneous periodic orbit
is Turing unstable after a sufficiently large rescaling of the
domain \cite[Theorem~3.5]{Yi-JDE-2021}.

We now give a local condition ensuring that the postcritical instability interval described in Proposition~\ref{pthirdprop-autonomous-periodic-modal-structure}
is nonempty. At the Hopf point, we have
\begin{equation}\label{pthird-Hopf-determinant}
	\det(\bar J-\alpha D)
	=
	d_1d_2\alpha^2
	-
	(d_2-d_1)\bar a\alpha
	+
	\Sigma.
\end{equation}
Define
\[
\bar \alpha_\pm
:=
\frac{
	(d_2-d_1)\bar a
	\pm
	\sqrt{
		\bigl[(d_2-d_1)\bar a\bigr]^2
		-
		4d_1d_2\Sigma
	}
}{
	2d_1d_2
}.
\]

\begin{proposition}\label{pthirdprop-hopf-finite-mode-persistence}
Suppose that
\begin{equation}\label{pthird-Hopf-finite-mode-condition}
		(d_2-d_1)\bar a>0,
		\qquad
		\bigl[(d_2-d_1)\bar a\bigr]^2
		>
		4d_1d_2\Sigma.
\end{equation}
Then for each fixed \(\alpha_0\in(\bar\alpha_-,\bar\alpha_+)\), we have $
\rho(T_{\varepsilon,\alpha_0})>1$ for all periodic orbits on the stable branch sufficiently close to the Hopf point. Consequently, on any fixed domain whose nonzero Neumann spectrum intersects \((\bar\alpha_-,\bar\alpha_+)\), all periodic orbits on the stable branch sufficiently close to the Hopf point are Turing unstable.
Such an intersection can be achieved by a suitable rescaling of any fixed domain.
\end{proposition}

\begin{proof}
We first consider the case \(d_2>d_1\). Condition
\eqref{pthird-Hopf-finite-mode-condition} implies \(\bar a>0\). Fix $\alpha_0\in(\bar\alpha_-,\bar\alpha_+)$. Since $\det(\bar J-\alpha_0D)<0$, the matrix \(\bar J-\alpha_0D\) has a positive eigenvalue. At the Hopf point, the corresponding Riccati equation is
\[
Z'
=
\bar q(1+Z^2)
-
\bigl[
2\bar a+\alpha_0(d_2-d_1)
\bigr]Z.
\]
It has two positive equilibria $Z_H^-<Z_H^+$. Since an equilibrium \(Z\) yields the eigenvalue \(\bar a-\alpha_0d_1-\bar q Z\)  of \(\bar J-\alpha_0D\), the smaller equilibrium \(Z_H^-\) corresponds to
its larger and positive eigenvalue. Hence $\bar a-\alpha_0d_1-\bar qZ_H^->0$. Choose $\overline Z\in(Z_H^-,Z_H^+)$ sufficiently close to \(Z_H^-\). Then \(\overline Z\) is a strict supersolution of the limiting Riccati equation and still satisfies
$\bar a-\alpha_0d_1-\bar q\,\overline Z>0$.

As \(\varepsilon\to\bar\varepsilon\) along the Hopf branch, the
coefficients of the transformed Riccati equation converge uniformly
to their values at the Hopf point. Therefore, both strict inequalities
above persist for all periodic orbits on the stable branch sufficiently
close to the Hopf point. The same constant \(\overline Z\) is therefore a
positive periodic supersolution with a positive growth
integral. Hence Remark~\ref{pthirdrem-supersolution} gives $\rho(T_{\varepsilon,\alpha_0})>1$. It also gives \(\alpha_0\ge\alpha_*^\varepsilon\), where \(\alpha_*^\varepsilon\) is the Riccati critical value for the chosen periodic orbit.
By Proposition~\ref{pthirdprop-autonomous-periodic-modal-structure}, the critical mode is strictly stable. Hence \(\alpha_0>\alpha_*^\varepsilon\), and \(\alpha_0\) belongs to the postcritical instability interval.

The case \(d_1>d_2\) follows by applying the same argument to the
rotated Riccati equation. The conclusions concerning the spatial
domain follow from the Neumann eigenvalue rescaling argument
in Lemma~\ref{pthirdlem-domain-scaling}.
\end{proof}

These two sets of sufficient conditions use different information.
The condition~\eqref{pthird-Yi-diagonal-Hopf-condition} involves
the nonlinear quantity
$\frac{\operatorname{Im}l_1}{\operatorname{Re}l_1}$ and detects instability through the phase multiplier at small positive modal parameters. By contrast, our condition~\eqref{pthird-Hopf-finite-mode-condition} depends
only on \(\bar J\) and \(D\), and ensures that the nearby stable periodic orbits have a nonempty postcritical instability interval. Neither of these two algebraic sufficient conditions implies the other in general. For instance, if \(\operatorname{Im}l_1=0\), then
\eqref{pthird-Yi-diagonal-Hopf-condition} cannot hold, whereas
\eqref{pthird-Hopf-finite-mode-condition} may still be satisfied.
Conversely, the condition~\eqref{pthird-Yi-diagonal-Hopf-condition} may hold even when the discriminant
of the polynomial in \eqref{pthird-Hopf-determinant} is nonpositive.

The distinction between these two criteria is also reflected in the
location of the instability regions in the modal parameter \(\alpha\).
When both \eqref{pthird-Yi-diagonal-Hopf-condition} and
\eqref{pthird-Hopf-finite-mode-condition} hold, fix a periodic
orbit on the stable branch sufficiently close to the Hopf point.
By Proposition~\ref{pthirdprop-autonomous-periodic-modal-structure},
its Riccati critical value is positive and has a neighborhood
of strictly stable modal parameters. The small-\(\alpha\) instability interval guaranteed 
the criterion~\eqref{pthird-Yi-diagonal-Hopf-condition} can be chosen to lie to the left of this neighborhood, whereas the nonempty postcritical instability interval lies to its right. Thus, the instability regions certified by the two conditions are separated by a nonempty interval of stable modal parameters.

\section{Applications}\label{pthirdsubsec-applications}

We now present two examples to illustrate our analytic results. The first concerns periodic orbits arising from a Hopf bifurcation, while the second concerns a nonconstant periodic solution in a periodic environment. In both examples, we impose homogeneous Neumann boundary conditions.
\subsection{The Sporns--Seelig model} \label{pthirdexp-sporns-seelig}
	
The Sporns--Seelig model describes an enzymatic reaction with substrate-induced enzyme synthesis \cite{Xiang-Wan-Lin-ERA-2022}. We consider the following choice of parameters in its reaction--diffusion version:
\begin{equation}\label{pthirdsys-sporns-seelig}
		\begin{cases}
			u_t=d_1\Delta u+1+2u-\displaystyle\frac{muv}{1+u},\\[2mm]
			v_t=d_2\Delta v+u-\displaystyle\frac{v}{1+u},
		\end{cases}
\end{equation}
where \(u\) and \(v\) denote the scaled substrate and enzyme concentrations, respectively, and \(m,d_1,d_2>0\). For \(s>0\), set
\[
	m=m(s):=\frac{1+2s}{s^2}.
\]
The kinetic system then has the positive equilibrium $E_s=\bigl(s,s(1+s)\bigr)$. By \cite[Theorem~2.1]{Xiang-Wan-Lin-ERA-2022}, a supercritical Hopf bifurcation occurs at \(s=1\). For \(s>1\) sufficiently close to \(1\), the kinetic system has a positive orbitally asymptotically stable periodic solution \((U_s(t),V_s(t))\), which converges uniformly to \((1,2)\) as \(s\to1^+\). Its least period \(\omega_s\) tends to \(\frac{2\pi}{\sqrt2}\) as \(s\to1^+\).

To match the sign convention \(b>0>c\), we apply the transformation \(P=\operatorname{diag}(1,-1)\), which reverses the signs of the off-diagonal entries while leaving the diagonal diffusion matrix unchanged. Then at the Hopf point, the transformed Jacobian is
\[
	\bar J=
	\begin{pmatrix}
		\frac12&\frac32\\[1mm]
		-\frac32&-\frac12
	\end{pmatrix},
	\qquad
	\bar a=\frac12,
	\qquad
	\Sigma=\det\bar J=2.
\]
Consequently, condition \eqref{pthird-Hopf-finite-mode-condition} becomes
\[
	d_2>d_1,
	\qquad
	(d_2-d_1)^2>32d_1d_2.
\]
This is equivalent to $\frac{d_2}{d_1}>17+12\sqrt2$. Thus, Proposition~\ref{pthirdprop-hopf-finite-mode-persistence} applies whenever the diffusion ratio satisfies this inequality.
	
In particular, take \(d_1=1\) and \(d_2=50\). Then
\[
	\det(\bar J-\alpha D)
	=
	50\alpha^2-\frac{49}{2}\alpha+2,
	\qquad
	\bar\alpha_\pm
	=
	\frac{49\pm\sqrt{801}}{200}.
\]
It is easy to verify that \(\alpha_0=1/5\) lies in \((\bar\alpha_-,\bar\alpha_+)\). Let \(T_{s,\alpha}\) denote the modal monodromy matrix along \((U_s,V_s)\). Proposition~
\ref{pthirdprop-hopf-finite-mode-persistence} gives $\rho(T_{s,1/5})>1$
for all \(s>1\) sufficiently close to \(1\).

To realize this unstable modal parameter by a genuinely two-dimensional
spatial mode, let $\Omega=(0,L)^2$ and $L=\pi\sqrt{10}$. The Neumann eigenfunctions of \(-\Delta\) on \(\Omega\) are
\[
\phi_{m,n}(x,y)
=
\cos\left(\frac{m\pi x}{L}\right)
\cos\left(\frac{n\pi y}{L}\right),
\qquad m,n=0,1,2,\ldots,
\]
with corresponding eigenvalues
\[
\lambda_{m,n}
=
\left(\frac{m\pi}{L}\right)^2
+
\left(\frac{n\pi}{L}\right)^2
=
\frac{m^2+n^2}{10}.
\]
In particular, we have $\lambda_{1,1}=\frac15$, and the corresponding eigenfunction is
\[
\phi_{1,1}(x,y)
=
\cos\left(\frac{x}{\sqrt{10}}\right)
\cos\left(\frac{y}{\sqrt{10}}\right).
\]
Hence the spatially homogeneous periodic solution \((U_s,V_s)\) is
unstable with respect to the \((1,1)\) Neumann mode. Since it is
orbitally asymptotically stable for the kinetic system, it is Turing
unstable for \eqref{pthirdsys-sporns-seelig} on \(\Omega\).

This example also illustrates the distinction from the sufficient condition in  \cite{Yi-JDE-2021}. With the normal form normalization in \cite[(2.12)]{Xiang-Wan-Lin-ERA-2022}, a direct computation at \(s=1\) gives
\[
	l_1=-\frac7{64}-\frac{\sqrt2}{64}i,
	\qquad
	\frac{\operatorname{Im}l_1}{\operatorname{Re}l_1}
	=\frac{\sqrt2}{7}.
\]
Therefore,
\[
	\sqrt{\Sigma}(d_1+d_2)
	+
	\frac{\operatorname{Im}l_1}{\operatorname{Re}l_1}
	(d_2-d_1)\bar a
	>0.
\]
Thus, condition~\eqref{pthird-Yi-diagonal-Hopf-condition} is not satisfied, and hence, the criterion in
 \cite{Yi-JDE-2021}
 is inconclusive for detecting Turing instability in this case.

\begin{figure}[!t]
	\centering
	\includegraphics[width=0.9\textwidth]{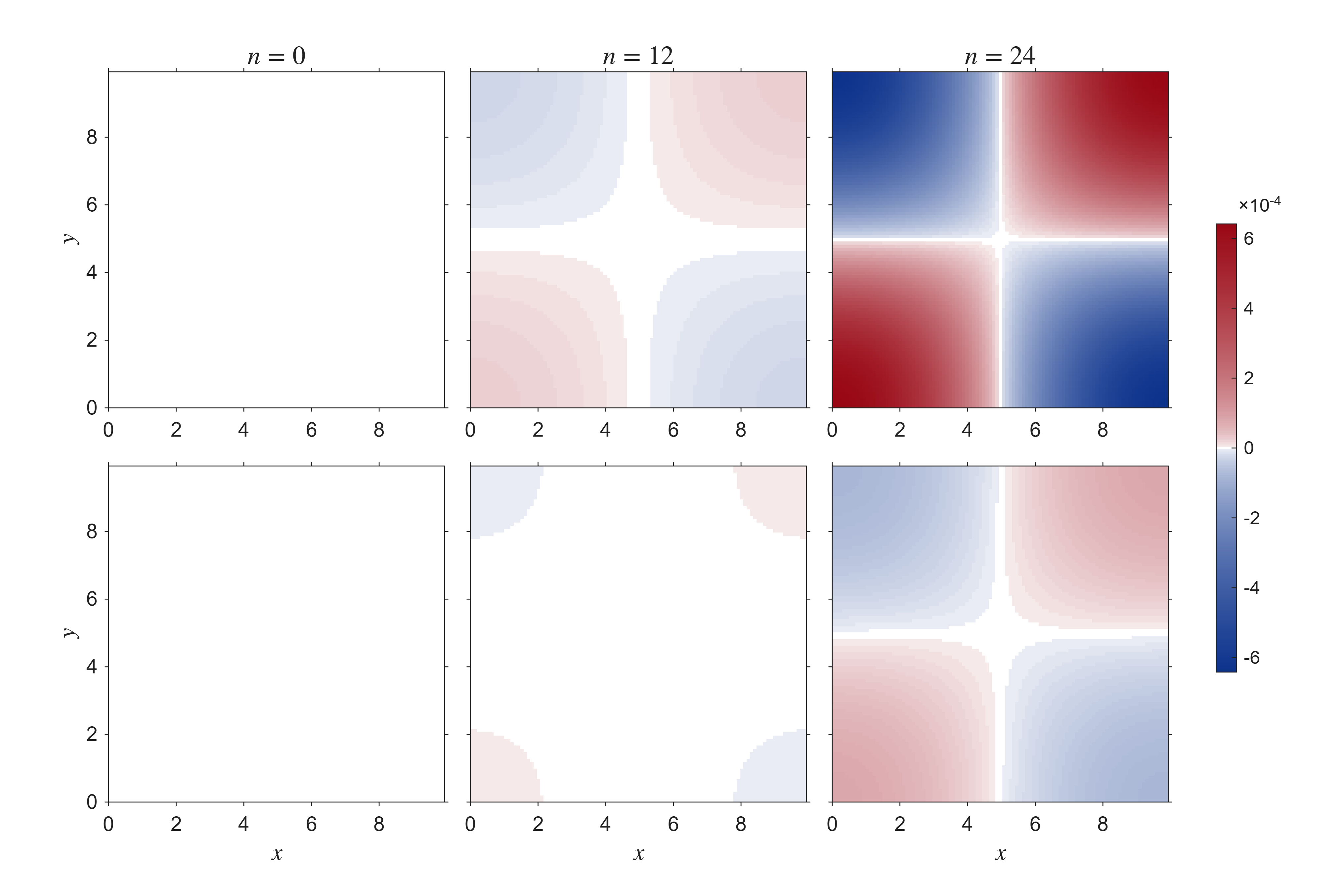}
	\caption{Numerical illustration of the Turing instability in
		the Sporns--Seelig model. The top row shows
		\(u(t,x,y)-U_s(t)\), and the bottom row shows
		\(v(t,x,y)-V_s(t)\), at
		\(t=0,\,12\omega_s,\) and $24\omega_s$,
		respectively. The initial perturbation is aligned with the unstable
		\((1,1)\) Neumann mode.}
	\label{pthirdfig-sporns-seelig}
\end{figure}

For a numerical illustration, we choose $s=1.02$. The corresponding kinetic periodic orbit has period $\omega_s\approx 4.49713840$. Numerically, the Floquet multipliers of the kinetic variational equation are $1$ and $0.87250963$. Thus the periodic orbit is orbitally
stable. For the spatial mode corresponding to $\alpha=\lambda_{1,1}=\frac15$, the modal monodromy matrix satisfies
$\rho(T_{s,1/5})\approx 1.30902176>1$. Hence a perturbation in the \((1,1)\) Neumann mode grows from one period to the next. To visualize this instability, we perturb the spatially homogeneous periodic solution in the direction of the dominant Floquet eigenvector associated with \(T_{s,1/5}\). Figure~\ref{pthirdfig-sporns-seelig} shows the deviations
from the spatially homogeneous periodic solution at the stroboscopic
times $t=0$, $12\omega_s$ and $24\omega_s$. The increasing amplitude of the deviations, while preserving the $(1,1)$ spatial structure, shows that this nonhomogeneous mode grows from period to period. Since the kinetic periodic orbit is stable, this growth of a nonzero spatial mode confirms the occurrence of Turing instability.

\subsection{The periodic Segel--Jackson model}\label{pthirdexp-segel-jackson}
	
We next consider the Segel--Jackson model in \cite[Section~5]{Sherratt-EJAM-1995}, with smooth periodic coefficients in place of the piecewise constant coefficients considered there:
\begin{equation}\label{pthirdsys-segel-jackson}
		\begin{cases}
			u_t=\Delta u+u+\displaystyle\frac12u^2-\mu_\varepsilon(t)uv,\\[2mm]
			v_t=d_2(t)\Delta v+uv-v^2,
		\end{cases}
\end{equation}
where \(u\) and \(v\) are the prey and predator densities, respectively, and
\[
	\mu_\varepsilon(t)=1+\varepsilon\cos(30t),
	\qquad
	d_2(t)=20+18\cos(30t),
	\qquad
	0\le\varepsilon<\frac12.
\]
The coefficients have period \(\omega=\pi/15\), and
\[
	d_2(t)-d_1(t)=19+18\cos(30t)\ge1.
\]
When \(\varepsilon=0\), the kinetic system has the positive equilibrium \((2,2)\), whose Jacobian is
\[
	J_0=
	\begin{pmatrix}
		1&-2\\
		2&-2
	\end{pmatrix}.
\]
It is easy to see that the equilibrium \(E_0=(2,2)\) is locally asymptotically stable since both eigenvalues of \(J_0\) have negative real parts.

Let \(\Pi_\varepsilon\) denote the Poincar\'e map of the kinetic system. At \(\varepsilon=0\), the equilibrium \(E_0\) is a fixed point of $\Pi_0$, whose derivative at \(E_0\) is \(e^{\omega J_0}\). Since the eigenvalues of \(J_0\) have negative real parts, \(1\notin\sigma(e^{\omega J_0})\). Hence the derivative of the corresponding fixed-point equation \(\Pi_0-Id\) is invertible at \(E_0\), where $Id$ is the identity map. By the implicit function theorem, there exists \(\varepsilon_0>0\) such that, for \(0<\varepsilon<\varepsilon_0\), the map \(\Pi_\varepsilon\) has a locally unique fixed point near \(E_0\). The corresponding solution, denoted by
\((U_\varepsilon,V_\varepsilon)\), is a positive \(\omega\)-periodic solution and converges uniformly to \(E_0\) on \([0,\omega]\) as \(\varepsilon\to0\). Moreover, by continuity of the Floquet multipliers, this periodic solution remains locally exponentially stable for sufficiently small \(\varepsilon>0\). This solution is nonconstant. Indeed, a positive constant solution would require \(U=V\) and \(1+\frac12U=\mu_\varepsilon(t)V\) for all \(t\), which is impossible when \(\varepsilon>0\).

For the present model, the system \eqref{pthird-linear-alpha-ODE} takes the form
\[
	\begin{pmatrix}
		\phi\\
		\psi
	\end{pmatrix}'
	=
	\begin{pmatrix}
		1+U_\varepsilon-\mu_\varepsilon V_\varepsilon-\alpha
		&
		-\mu_\varepsilon U_\varepsilon\\[1mm]
		V_\varepsilon
		&
		U_\varepsilon-2V_\varepsilon-\alpha d_2(t)
	\end{pmatrix}
	\begin{pmatrix}
		\phi\\
		\psi
	\end{pmatrix}.
\]
Since $(U_\varepsilon,V_\varepsilon)$ satisfies the kinetic system, we obtain
\[
	\frac{U_\varepsilon'}{U_\varepsilon}
	=
	1+\frac12U_\varepsilon-\mu_\varepsilon V_\varepsilon,
	\qquad
	\frac{V_\varepsilon'}{V_\varepsilon}
	=
	U_\varepsilon-V_\varepsilon.
\]
Introduce $\phi=U_\varepsilon(t)\xi$ and $\psi=-V_\varepsilon(t)\eta$. Then the modal equation becomes
	
\begin{equation}\label{pthirdsys-segel-jackson-modal}
		\begin{cases}
			\xi'
			=
			\left(\displaystyle\frac12U_\varepsilon-\alpha\right)\xi
			+\mu_\varepsilon V_\varepsilon\eta,\\[2mm]
			\eta'
			=
			-U_\varepsilon\xi
			+\left(-V_\varepsilon-\alpha d_2(t)\right)\eta.
		\end{cases}
\end{equation}
This periodic invertible transformation preserves the Floquet multipliers. For the coefficients of \eqref{pthirdsys-segel-jackson-modal}, we have
\[
	r_\varepsilon
	=
	\sqrt{\frac{\mu_\varepsilon V_\varepsilon}{U_\varepsilon}},
	\qquad
	q_\varepsilon
	=
	\sqrt{\mu_\varepsilon U_\varepsilon V_\varepsilon},
	\qquad
	\gamma_\varepsilon
	=
	\frac{r_\varepsilon'}{r_\varepsilon}.
\]
The corresponding Riccati equation is
\[
	Z'
	=
	F_{\varepsilon,\alpha}(t,Z)
	:=
	q_\varepsilon(1+Z^2)
	-
	\left[
	\frac12U_\varepsilon+V_\varepsilon
	+\alpha\bigl(d_2(t)-1\bigr)
	-\gamma_\varepsilon
	\right]Z,
\]
and the associated growth integral is
\[
	I_{\varepsilon,\alpha}[Z]
	=
	\int_0^\omega
	\left(
	\frac12U_\varepsilon-\alpha-q_\varepsilon Z
	\right)\,\mathrm dt.
\]

Choose
\[
	\alpha_0=\frac3{10},
	\qquad
	\vartheta(t)
	=
	\frac3{10}
	\exp\left[-\frac9{50}\sin(30t)\right],
	\qquad
	\overline Z_\varepsilon(t)
	=
	r_\varepsilon(t)\vartheta(t).
\]
Since
\[
	\vartheta'(t)
	=
	-\frac{27}{5}\cos(30t)\vartheta(t),
\]
a direct calculation gives
\[
	F_{\varepsilon,\alpha_0}\bigl(t,\overline Z_\varepsilon(t)\bigr)
	-\overline Z_\varepsilon'(t)=r_\varepsilon(t)
	\left[
	U_\varepsilon
	+\mu_\varepsilon V_\varepsilon\vartheta^2
	-
	\left(
	\frac12U_\varepsilon+V_\varepsilon+\frac{57}{10}
	\right)\vartheta
	\right].
\]
At \(\varepsilon=0\), the expression in brackets is $2(1+\vartheta^2)-\frac{87}{10}\vartheta$. Since \(1/4<\vartheta(t)<3/8\), we have
\[
	2(1+\vartheta(t)^2)-\frac{87}{10}\vartheta(t)
	\le-\frac1{20}<0
	\qquad\text{for all }t\in[0,\omega].
\]
Thus, the supersolution inequality is strict at \(\varepsilon=0\).
	
Furthermore, \(q_\varepsilon r_\varepsilon=\mu_\varepsilon V_\varepsilon\), and hence
\[
	I_{\varepsilon,\alpha_0}[\overline Z_\varepsilon]
	=
	\int_0^\omega
	\left[
	\frac12U_\varepsilon-\frac3{10}
	-\mu_\varepsilon V_\varepsilon\vartheta
	\right]\,\mathrm dt.
\]
At \(\varepsilon=0\), convexity of the exponential function and the zero mean of \(\sin(30t)\) yield
\[
	I_{0,\alpha_0}[\vartheta]
	=\int_0^\omega
	\left(\frac7{10}-2\vartheta(t)\right)\,\mathrm dt
	\ge
	\left[
	\frac7{10}-\frac35\cosh\left(\frac9{50}\right)
	\right]\omega
	> 
	\frac9{100}\,\omega
	>0.
\]
Both strict inequalities persist for sufficiently small \(\varepsilon>0\). After reducing \(\varepsilon_0\) if necessary, \(\overline Z_\varepsilon\) is therefore a positive periodic supersolution satisfying $I_{\varepsilon,\alpha_0}[\overline Z_\varepsilon]>0$ for $0<\varepsilon<\varepsilon_0$. Remark~\ref{pthirdrem-supersolution} gives \(\rho(T_{\varepsilon,\alpha_0})>1\).

To realize this unstable modal parameter by a genuinely two-dimensional
spatial mode, let $\Omega=(0,L)^2$ and $L=\pi\sqrt{\frac{20}{3}}$. The Neumann eigenfunctions of \(-\Delta\) on \(\Omega\) are
\[
\phi_{m,n}(x,y)
=
\cos\left(\frac{m\pi x}{L}\right)
\cos\left(\frac{n\pi y}{L}\right),
\qquad m,n=0,1,2,\ldots,
\]
with corresponding eigenvalues
\[
\lambda_{m,n}
=
\left(\frac{m\pi}{L}\right)^2
+
\left(\frac{n\pi}{L}\right)^2
=
\frac{3(m^2+n^2)}{20}.
\]
In particular, we have $\lambda_{1,1}=\frac3{10}=\alpha_0$, and the corresponding eigenfunction is
\[
\phi_{1,1}(x,y)
=
\cos\left(\sqrt{\frac3{20}}\,x\right)
\cos\left(\sqrt{\frac3{20}}\,y\right).
\]
Hence \((U_\varepsilon,V_\varepsilon)\) is unstable with respect to
the \((1,1)\) Neumann mode. Since it is locally exponentially stable
for the kinetic system, Theorem~\ref{pthirdthm-refined-turing-instability}
implies that the nonconstant periodic solution
\((U_\varepsilon,V_\varepsilon)\) is Turing unstable for
\eqref{pthirdsys-segel-jackson} whenever
\(0<\varepsilon<\varepsilon_0\).

\begin{figure}[!htbp]
	\centering
	\includegraphics[width=0.9\textwidth]{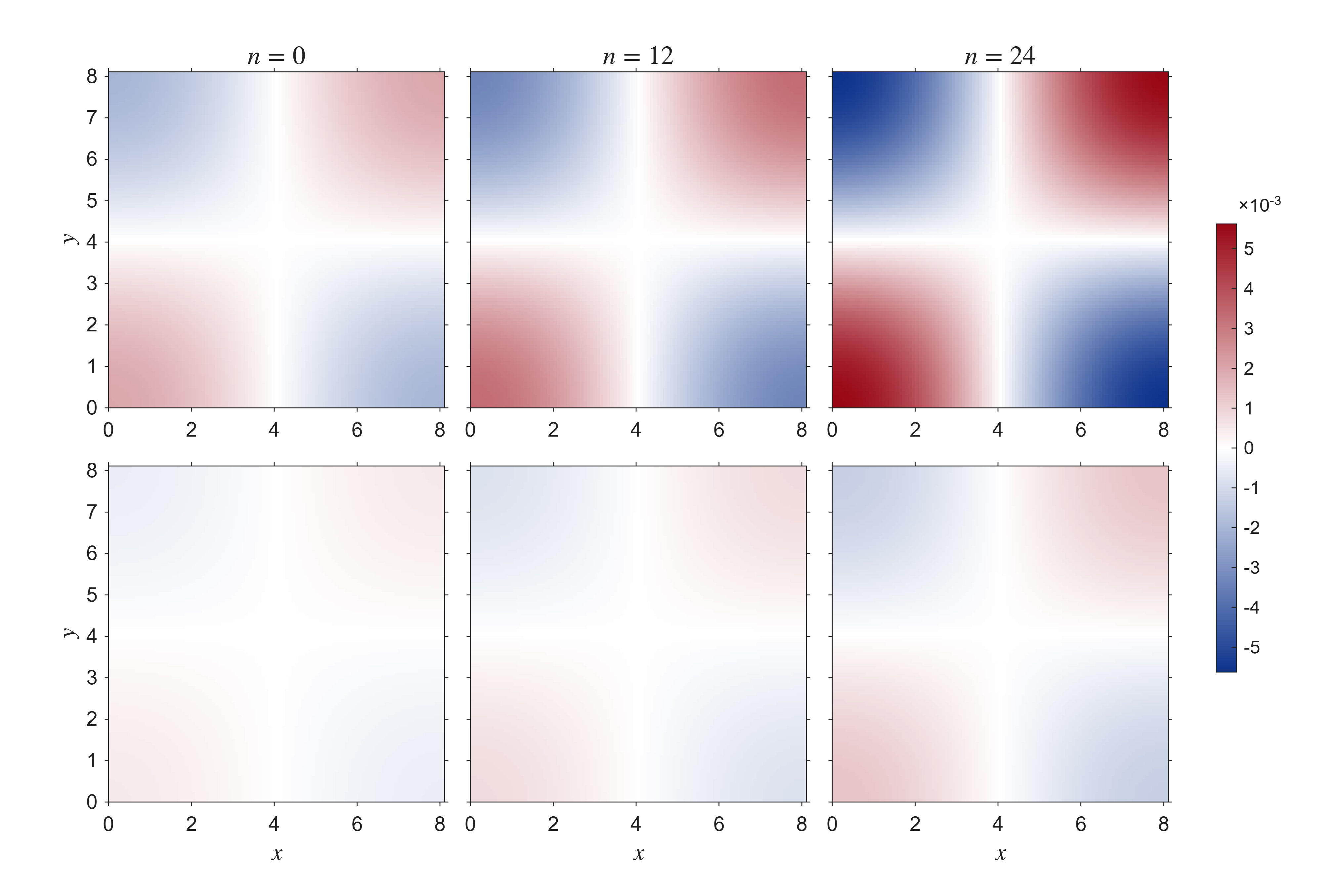}
	\caption{Numerical illustration of the Turing instability in the
		periodic Segel--Jackson model.
		The top row shows
		\(u(t,x,y)-U_\varepsilon(t)\), and the bottom row shows
		\(v(t,x,y)-V_\varepsilon(t)\), at
		\(t=0,\,12\omega\), and \(24\omega\),
		respectively. The initial perturbation is aligned with the unstable
		\((1,1)\) Neumann mode.}
	\label{pthirdfig-segel-jackson}
\end{figure}

We next illustrate the instability in the time-periodic setting by taking
$\varepsilon=0.02$. The corresponding kinetic periodic solution has period
$\omega=\pi/15\approx 0.20943951$. Numerically, the Floquet multipliers of
the kinetic variational equation are $0.86623174\pm0.24633631i$, with spectral radius
$0.90057704<1$. Thus the periodic solution is locally exponentially stable.
For the spatial mode corresponding to
$\alpha=\lambda_{1,1}=\frac3{10}$, the modal monodromy matrix satisfies
$\rho(T_{\varepsilon,3/10})\approx1.04394789>1$.
Hence a perturbation in the \((1,1)\) Neumann mode grows from one period
to the next. As in the previous example, the instability is visualized by perturbing the spatially homogeneous periodic solution in the direction of the dominant Floquet eigenvector associated with \(T_{\varepsilon,3/10}\).
Figure~\ref{pthirdfig-segel-jackson} shows the deviations from the spatially
homogeneous periodic solution at the stroboscopic times
$t=0$, $12\omega$ and $24\omega$.

Finally, we verify the comparison observation in
Remark~\ref{pthirdrem-Van-Gorder}. Let \(t_*=\omega/2=\pi/30\).
Since \(d_2(t_*)=2\) and
\(\mu_\varepsilon'(t_*)=d_2'(t_*)=0\), at
\(\varepsilon=0\) and \(\alpha_0=3/10\) we have
\[
\det\bigl(A_0(t_*)-\alpha_0D(t_*)\bigr)
=\frac{109}{50}>0,
\]
whereas
\[
b_0(t_*)\frac{\mathrm d}{\mathrm dt}
\left(\frac{a_0-\alpha_0}{b_0}\right)\Big|_{t=t_*}
=
c_0(t_*)\frac{\mathrm d}{\mathrm dt}
\left(\frac{e_0-\alpha_0d_2}{c_0}\right)\Big|_{t=t_*}
=0.
\]
By continuity, both pointwise conditions in
Remark~\ref{pthirdrem-Van-Gorder} therefore fail at \(t=\omega/2\)
for all sufficiently small \(\varepsilon>0\).

\section*{Statements and Declarations}
\noindent\textbf{Funding.}
Liang's research was supported in part by the
China Scholarship Council (202407970004), and Zhao's research was supported in part by the Natural Sciences and Engineering Research Council of Canada (RGPIN-2025-04963).

\medskip
\noindent\textbf{Data availability.}
No external datasets were used in this study.

\medskip
\noindent\textbf{Conflict of interest.}
On behalf of all authors, the corresponding author states that there is no conflict of interest.

\end{document}